\RequirePackage{ifpdf}
\ifpdf 
\documentclass[pdftex]{sigma}
\else
\documentclass{sigma}
\fi

\let\ds=\displaystyle
\def\coun{\varepsilon}

\def\<{\langle}
\def\>{\rangle}

\newcommand{\ZZ}{{\mathbb Z}}

\def\r#1{(\ref{#1})}
\let\rf=\r
\def\ot{\otimes}

\def\sk#1{\left(#1\right)}

\def\U{\overline{U}}

\def\Pfpm{{P}^\pm}
\def\Pfp{{P}^+}
\def\Pfm{{P}^-}

\def\Uqgln{U_q(\widehat{\mathfrak{gl}}_N)}

\def\Uqgl#1{U_q(\widehat{\mathfrak{gl}}_{#1})}

\def\Uqsl2{U_q(\widehat{\mathfrak{sl}}_2)}

\def\Uqbp{U_q(\mathfrak{b}^+)}
\def\ff{{F}}

\def\Uqbp{U_q(\mathfrak{b}^+)}

\def\si{\sigma}

\def\LL{{\rm L}}

\def\F{{\mathcal{F}}}

\def\tSym{\overline{\rm Sym}}

\def\Uqqq{U_q(\mathfrak{gl}_N)}

\def\lll{l}
\def\rr{r}
\def\ss{s}
\def\mm{m}

\def\seg#1#2{[#2,#1]}
\def\segg#1{[#1]}
\def\meg#1#2{[#1,#2]}
\def\megg#1{[#1]}

\def\Ff{f}

\def\bFF{\mathbf{F}}

\def\Fgs{\mathcal{F}}
\def\Ags{\mathcal{A}}
\def\Bgs{\mathcal{B}}
\def\Cgs{\mathcal{C}}
\def\Sgs{\mathcal{S}}
\def\tSgs{\tilde{\mathcal{S}}}

\def\heit{{\sf h}}
\def\tchi{{\bar\chi}}
\def\wt{\bar n}
\def\Ccal{\mathcal{E}}
\def\Dcal{\mathcal{D}}

\numberwithin{equation}{section}

\begin{document}

\allowdisplaybreaks

\renewcommand{\thefootnote}{$\star$}

\renewcommand{\PaperNumber}{081}

\FirstPageHeading

\ShortArticleName{Generating Series for Nested Bethe Vectors}

\ArticleName{Generating Series for Nested Bethe Vectors\footnote{This paper is a
contribution to the Special Issue on Kac--Moody Algebras and Applications. The
full collection is available at
\href{http://www.emis.de/journals/SIGMA/Kac-Moody_algebras.html}{http://www.emis.de/journals/SIGMA/Kac-Moody{\_}algebras.html}}}

\Author{Sergey KHOROSHKIN~$^\dag$ and Stanislav PAKULIAK~$^{\ddag\dag}$}

\AuthorNameForHeading{S. Khoroshkin  and S. Pakuliak}

\Address{$^\dag$~Institute of Theoretical {\rm \&} Experimental Physics, 117259 Moscow, Russia} 
\EmailD{\href{mailto:khor@itep.ru}{khor@itep.ru}}

\Address{$^\ddag$~Laboratory of Theoretical Physics, JINR, 141980 Dubna, Moscow Region, Russia}
\EmailD{\href{mailto:pakuliak@theor.jinr.ru}{pakuliak@theor.jinr.ru}}

\ArticleDates{Received September 14, 2008; Published online November 24, 2008}

\vspace{-1mm}

\Abstract{We reformulate  nested relations between of\/f-shell $U_q(\widehat{\mathfrak{gl}}_N)$ Bethe
vectors as a~certain equation on generating series of
strings of the composed  $U_q(\widehat{\mathfrak{gl}}_N)$ currents.  Using inversion of the
generating series we f\/ind a new type of hierarchical relations
between univer\-sal of\/f-shell Bethe vectors, useful for a derivation
of Bethe equation. As an example of application,  we use these
relations for a derivation of  analytical Bethe ansatz equations
[Arnaudon D. et al.,
{\it  Ann. Henri  Poincar\'e} {\bf 7} (2006), 1217--1268, \href{http://arxiv.org/abs/math-ph/0512037}{math-ph/0512037}] for the para\-meters of  universal Bethe vectors of the algebra $U_q(\widehat{\mathfrak{gl}}_2)$.}

\vspace{-1mm}


\vspace{-1mm}

\Classification{17B37; 81R50}

\renewcommand{\thefootnote}{\arabic{footnote}}
\setcounter{footnote}{0}

\vspace{-4mm}

\section{Introduction}

Hierarchical (nested) Bethe ansatz (NBA) was designed \cite{KR83} to solve quantum integrable mo\-dels with
$\mathfrak{gl}_N$ symmetries.
The cornerstone of NBA is a procedure which relates  Bethe vectors for the
model with $\mathfrak{gl}_N$ symmetry to  analogous objects with
 $\mathfrak{gl}_{N-1}$ symmetry.
 This hierarchical procedure is implicit, it allows to obtain
 Bethe equations for the parameters of the Bethe vectors while the explicit construction of these vectors
 itself remains rather non-trivial problem.
Authors of the papers \cite{VT,VT1} proposed a closed expression for
 of\/f-shell Bethe vectors as matrix elements of the monodromy
operator built of products of  fundamental $\LL$-operators. However
a calculation of this expression in every particular representation
 is still a  nontrivial problem. Such a calculations was done in \cite{VT2} on the level of the
evaluation homomorphism of $\Uqgln\to\Uqqq$.

The construction of \cite{VT,VT1} yields the of\/f-shell Bethe vectors in terms of  matrix elements
of the monodromy matrix which satisf\/ies the corresponding quantum Yang--Baxter equation
 and generate a Borel subalgebra of quantum af\/f\/ine algebra $\Uqgln$ \cite{RS} for trigonometric $R$-matrix
 or the Yangian $Y(\mathfrak{gl}_N)$   for the rational $R$-matrix.
Those quantum af\/f\/ine algebras as well as doubles of Yangians possess another ``new'' realization
introduced in \cite{D88}. In this realization the corresponding algebra is described in terms
of  generating series (currents) and an isomorphism between dif\/ferent realizations of these
inf\/inite-dimensional algebras was observed in \cite{DF}.
Using this isomorphism one may try to look for the expressions for the universal of\/f-shell Bethe vectors
in terms of the modes of the currents. This program was realized in
\cite{KP,KPT},  where  explicit formulas for the of\/f-shell Bethe vectors
in terms of the currents were found.
A signif\/icant part of this  approach to the construction of  Bethe vectors
is a method of projection introduced in \cite{ER} and developed in the recent paper \cite{EKhP}.
This method operates with projections of Borel subalgebra to its intersections with  Borel
subalgebras of a dif\/ferent type.
It was proved in \cite{EKhP} that  universal of\/f-shell Bethe vectors can be identif\/ied
with the projections of  products of the Drinfeld currents to the intersections of
Borel subalgebras of dif\/ferent
types. In the latter paper it was checked in a rather general setting  that the Bethe vectors
obtained from the projections of the currents satisfy the same comultiplication rule as
the Bethe vectors constructed in terms of the fundamental $\LL$-operators.
This approach was used in \cite{KhP-GLN,OPS} for further generalization of the results obtained
in \cite{VT2}.

A deduction of  general expressions for of\/f-shell Bethe vectors  \cite{KhP-GLN,OPS} is based on
hierarchical relations between projections of  products of the currents (see Proposition~4.2 in
\cite{KhP-GLN}). However these relations do not help in the investigation of the action of
integrals of motions on of\/f-shell Bethe vectors. The situation is quite dif\/ferent in classical
approach of nested Bethe ansatz. The corresponding hierarchical relations allow to compute the action
of integrals of motion  and derive the Bethe equations \cite{KR83,MTV}, but can be hardly used in
the computation of the explicit expressions of Bethe vectors. This signif\/ies an existence of two
types of hierarchical relations for of\/f-shell Bethe vectors and thus of two dif\/ferent presentations of them.
The goal of this paper is to observe these two type of the hierarchical relations within the
approach of the method of projections.

The  new important  objects which appear in the application of the method of projections to the investigation
of the of\/f-shell Bethe vectors are so called {\it strings} and their projections. Strings are special
ordered products of composed current introduced in \cite{DKh}. Our hierarchical relations express
$\Uqgln$ Bethe vectors via  products of  special strings and $\Uqgl{N-1}$ Bethe vectors. On the other
hand the basic point of application of the method of projection is the ordered decomposition of the product
of the currents. The factors  of this decomposition are   projections of the products of the currents
\r{dec-ff1m}. The crucial observation is that this decomposition and hierarchical relations
for the opposite projections of the products of currents have a similar structure. Combinations of relations of
these two dif\/ferent type
can be used to obtain a new type hierarchical relations for of\/f-shell Bethe vectors.

\looseness=-1
In order to solve the latter problem we collect all of\/f-shell Bethe vectors into multi-variable
generating series. We also introduce the generating series for the products of the currents and
for  the strings and their projections. We rewrite the hierarchical relations
for $\Uqgln$ Bethe vectors  as a simple relation
on the product of the generating series of the projections of the strings and
the generating series of the $\Uqgl{N-1}$ Bethe vectors.
A similar construction can be performed for of\/f-shell Bethe vectors related to opposite Borel
subalgebra.
However the product entering into these
relations is not usual. It contains a $q$-symmetrization with a special functional weight.
This $\star$-product is associative and the generating series are invertible with respect
to this product. Finally the new type of the hierarchical relations reduces to inversion
of the generating series of opposite projections of the strings. This inversion is ef\/fectively performed.
To do this we introduce some combinatorial language of tableaux f\/illed by the Bethe parameters
(see Section~\ref{combi}).

Applications of the new type hierarchical relations to the investigation of the properties of
quantum integrals of motion are given in \cite{FKPR}. Here we demonstrate how they work in the
simplest case of the universal $\Uqgl{2}$ Bethe vectors. As a result we get universal Bethe
equations of the analytical Bethe ansatz \cite{ACDFR}. In contrast to the usual Bethe equations
these equations refer to Cartan currents instead of the eigenvalues of the diagonal elements
of monodromy matrix on highest weight vectors. In the Appendix we collect the basic def\/ining relations
for the $\Uqgln$ composed currents.

\section[Generating series  and $\star$-product]{Generating series  and $\boldsymbol{\star}$-product}

\subsection[A $q$-symmetrization]{A $\boldsymbol{q}$-symmetrization}

Let $\bar t=\{t_1,\ldots,t_n\}$ be a set of formal variables. Let $G(\bar t)$
be a Laurent series taking values in $\Uqgln$.
Consider the permutation group $S_n$ and its action on the formal series of $n$ variables
$\bar t=\{t_1,\ldots,t_n\}$
def\/ined for the elementary transpositions $\sigma_{i,i+1}$ as follows
\begin{align}\label{tr-fac}
\pi(\sigma_{i,i+1})G(t_1, \dots,t_i, t_{i+1},\dots, t_n) =
\frac{q^{-1}-q t_i/t_{i+1}}{q-q^{-1} t_i/t_{i+1}}  G(t_1, \dots,t_{i+1}, t_i, \dots, t_n),
\end{align}
where the rational series $\frac{1}{1-x}$ is understood as a series $\sum_{n\geq 0}x^n$
and $q$ is a deformation parameter of $\Uqgln$.
 Summing the action over  the group of permutations
we obtain the operator $\tSym_{\,\bar t}=\sum_{\sigma\in S_n}\pi(\sigma)$ acting as
follows:
\begin{align}\label{qss}
\tSym_{\, \bar t}\  G(\bar t) =\frac{1}{n!}
\sum\limits_{\si \in S_n}\prod\limits_{\substack{\ell<\ell'\\ \si(\ell)>\si(\ell')}}
\frac{q^{-1}-q t_{\si(\ell')}/t_{\si(\ell)}}
{q-q^{-1} t_{\si(\ell')}/t_{\si(\ell)}}   G(^\sigma t).
\end{align}
The product is taken over all pairs $(\ell, \ell')$, such that conditions $\ell
< \ell'$ and $\si(\ell) > \si(\ell')$ are satisf\/ied simultaneously.

We call operator  $\tSym_t$  a {\em $q$-symmetrization}.
One can check that the operation given by  \r{qss} is a projector
\begin{gather}\label{sym*}
\tSym_{\, \bar t}\, \tSym_{\, \bar t}\, (\cdot)= \tSym_{\, \bar t}\, (\cdot).
\end{gather}

Fix any positive integer $N>1$. Let $\bar\lll=\{\lll_1,\ldots,\lll_{N-1}\}$ and
$\bar\rr=\{\rr_1,\ldots,\rr_{N-1}\}$ be the sets of non-negative integers
satisfying a set of inequalities
\begin{gather}\label{set23}
\lll_a\leq\rr_a,\qquad a=1,\ldots,N-1.
\end{gather}
Denote by $\meg{\bar\lll}{\bar\rr}$  a set of segments which contain  positive
integers $\{\lll_a+1,\lll_a+2,\ldots,\rr_a-1,\rr_a\}$
including $\rr_a$ and excluding $\lll_a$. The length of each segment is equal
to $\rr_a-\lll_a$.

For a given set $\meg{\bar\lll}{\bar\rr}$ of segments  we denote by
$\bar t_{\meg{\bar\lll}{\bar\rr}}$  the sets of variables
\begin{gather}\label{t-coll}
\bar t_{\seg{\bar\rr}{\bar\lll}} =
\big\{t^{1}_{\lll_{1}+1},\ldots,t^{1}_{\rr_{1}};
t^{2}_{\lll_{2}+1},\ldots,t^{2}_{\rr_{2}};\ldots;
t^{N-1}_{\lll_{N-1}+1},\ldots,t^{N-1}_{\rr_{N-1}} \big\}.
\end{gather}
For any $a=1,\ldots,N-1\,$ we denote the sets of variables corresponding to the  segments
$\meg{\lll_a}{\rr_a}=\{\lll_a+1,\lll_a+2,\ldots,\rr_a\}$
as $\bar{t}^a_{\meg{{\lll_a}}{\rr_a}}\
=\{t^{a}_{\lll_{a}+1},\ldots,t^{a}_{\rr_{a}} \}$.
All the variables in $\bar{t}^a_{\meg{{\lll_a}}{\rr_a}}$ have  the type $a$. For the segments
$\meg{\lll_a}{\rr_a}=\meg{0}{n_a}$ we use the shorten  notations
$\bar{t}_{\meg{\bar 0}{\bar n}}\equiv\bar{t}_{\megg{\bar n}}$ and
$\bar{t}^a_{\meg{0}{n_a}}\equiv\bar{t}^a_{\megg{n_a}}$. We also name for a short collections
$\meg{\bar\lll}{\bar\rr}$ of segments as a {\it segment}.

Denote by $S_{\bar\lll,\bar\rr} =
S_{\lll_{1},\rr_{1}}\times \cdots \times S_{\lll_{N-1},\rr_{N-1}}$
 a direct product of the groups $S_{\lll_{a},\rr_{a}}$ permuting
integers  $\lll_{a}+1,\ldots, \rr_{a}$.
Let $G(\bar t_{\seg{\bar\rr}{\bar\lll}})$ be a series
depending on the ratios $t^a_i/t^b_j$ for $a<b$ and  $t^a_i/t^a_j$ for $i<j$.
The $q$-sym\-me\-tri\-za\-tion over the whole set of variables
$\bar t_{\seg{\bar\rr}{\bar\lll}}$  of the series $G(\bar t_{\seg{\bar\rr}{\bar\lll}})$
 is def\/ined by the formula
\begin{gather}\label{qsr}
\tSym_{\, \bar t_{\seg{\bar\rr}{\bar\lll}}} \, G(\bar
t_{\seg{\bar\rr}{\bar\lll}})= \sum_{\si\in
S_{\bar\lll,\bar\rr}}\prod_{1\leq a\leq N-1}\frac{1}{(\rr_a-\lll_a)!}
\prod_{\substack{\ell<\ell'\\ \si^a(\ell)>\si^a(\ell')}}
\frac{q^{-1}-q^{}\,t^a_{\si^a(\ell')}/t^a_{\si^a(\ell)}}
{q^{}-q^{-1} t^a_{\si^a(\ell')}/t^a_{\si^a(\ell)}}  G(^\si \bar
t_{\seg{\bar\rr}{\bar\lll}}),\!\!
\end{gather}
where  the set $^\si \bar t_{\seg{\bar\rr}{\bar\lll}}$
is def\/ined  as
\begin{gather}\label{sigmat}
^\si \bar t_{\seg{\bar\rr}{\bar\lll}} =
\big\{t^{1}_{\si^1(\lll_{1}+1)},\ldots,t^{1}_{\si^1(\rr_{1})};
t^{2}_{\si^2(\lll_{2}+1)},\ldots,t^{2}_{\si^2(\rr_{2})};\ldots;
t^{N-1}_{\si^{N-1}(\lll_{N-1}+1)},\ldots,t^{N-1}_{\si^{N-1}(\rr_{N-1})} \big\}.\!\!\!
\end{gather}

We say that the series $G(\bar t_{\seg{\bar\rr}{\bar\lll}})$ is {\it $q$-symmetric}, if
it is invariant under the action $\pi$ of each group
$S_{\lll_{a},\rr_{a}}$ with respect to the permutations of the variables
$t^a_{\lll_a+1},\ldots,t_{\rr_a}$ for $a=1,\ldots,N-1$:
\begin{gather}\label{exa3}
\tSym_{\, \bar t_{\seg{\bar\rr}{\bar\lll}}}
G(\bar t_{\seg{\bar\rr}{\bar\lll}})=
G(\bar t_{\seg{\bar\rr}{\bar\lll}}).
\end{gather}
Due to \r{sym*} the $q$-symmetrization $G(\bar t_{\seg{\bar\rr}{\bar\lll}})=
\tSym_{\,\bar t_{\seg{\bar\rr}{\bar\lll}}} Q(\bar t_{\seg{\bar\rr}{\bar\lll}})$
of any series $ Q(\bar t_{\seg{\bar\rr}{\bar\lll}})$ is a $q$-sym\-met\-ric series.

The rational function
\begin{gather}\label{ex-qsym}
\prod_{a=1}^{N-1} \prod_{i<j} \frac{1-t^a_i/t^a_j}{q-q^{-1}\,t^a_i/t^a_j}
\end{gather}
which is understood as a series with respect to $t^a_i/t^a_j$ is an example of a $q$-symmetric series.

\subsection{Generating series}

Let $u_i$, $i=1,\ldots,N-1$ be formal parameters. We denote the set of these parameters as
$\bar u=\{u_1,\ldots,u_{N-1}\}$.
Def\/ine a generating series
\begin{gather}\label{gs5}
\Ags(\bar u)=1+\sum_{n_1,\ldots,n_{N-1}\geq0\atop {\bar n}\neq{\bar 0}} \Ags(\bar t_{\segg{\bar n}})\
u_1^{n_1} u_2^{n_2}\cdots u_{N-1}^{n_{N-1}},
\end{gather}
where the coef\/f\/icients $\Ags(\bar t_{\segg{\bar n}})$ are arbitrary $q$-symmetric
series ($\tSym_{\, \bar t_{\segg{\bar n}}} \sk{\Ags(\bar t_{\segg{\bar n}})}=
\Ags(\bar t_{\segg{\bar n}}) $) of the
formal variables $\bar t_{\bar n}$ numbered by the multi-index $\bar n=\{n_1,\ldots,n_{N-1}\}$:
\begin{gather}\label{set111}
\bar{t}_{\segg{\bar{n}}} = \left\{ t^{1}_{1},\ldots,t^{1}_{n_{1}};
t^{2}_{1},\ldots, t^{2}_{n_{2}}; \ldots ;
t^{N-2}_{1},\ldots, t^{N-2}_{n_{N-2}};
t_1^{N-1},\ldots,t^{N-1}_{n_{N-1}}\right\}.
\end{gather}
We call  generating series of this type  {\it $q$-symmetric generating series}.
Note that the multi-index~$\bar n$ of the coef\/f\/icients of a $q$-symmetric generating series
is uniquely def\/ined by the set of formal variables $\bar t_{\segg{\bar n}}$, thus the coef\/f\/icients
$\Ags(\bar t_{\segg{\bar n}})$ are used without any additional index. However once it will be convenient
for us to use instead the notation $\Ags_{\bar n}\equiv\Ags(\bar t_{\segg{\bar n}})$
(see proof of Proposition~\ref{inversion}).

For two generating
series $\Ags(\bar u)$ and  $\Bgs(\bar u)$ we def\/ine  $\star$-product as a generating series
\begin{gather}\label{gs6}
\Cgs(\bar u)=\Ags(\bar u)\,\star\,\Bgs(\bar u)=
\sum_{n_1,\ldots,n_{N-1}\geq0} \Cgs(\bar t_{\segg{\bar n}})\,
u_1^{n_1}\,u_2^{n_2}\cdots u_{N-1}^{n_{N-1}}
\end{gather}
 with coef\/f\/icients
\begin{gather}\label{gs7}
 \Cgs(\bar t_{\segg{\bar n}})=
\sum_{0\leq \ss_{N-1}\leq n_{N-1}}\cdots \sum_{0\leq \ss_1\leq n_1}
   \tSym_{\, \bar t_{\segg{\bar n}}}
\left(Z_{\bar\ss}({\bar t}_{\segg{\bar n}})
  \Ags(\bar t_{\seg{\bar n}{\bar\ss}})\cdot \Bgs(\bar t_{\seg{\bar\ss}{\bar0}})\right),
\end{gather}
where
\begin{gather}\label{anm}
\Ags(\bar t_{\seg{\bar n}{\bar\ss}}) =
\Ags\big(t^{1}_{\ss_{1}+1},\ldots,t^{1}_{n_{1}};
t^{2}_{\ss_{2}+1},\ldots,t^{2}_{n_{2}};\ldots;
t^{N-1}_{\ss_{N-1}+1},\ldots,t^{N-1}_{n_{N-1}}\big)
\end{gather}
and
\begin{gather}\label{Zserm}
Z_{\bar\ss}({\bar t}_{\segg{\bar n}})=\prod_{a=1}^{N-2}
\prod_{\substack{ \ss_a < \ell\leq n_a \\
  0< \ell' \leq \ss_{a+1}}} \frac{q-q^{-1}
t^{a}_{\ell}/ t^{a+1}_{\ell'}}{1-t^{a}_{\ell}/t^{a+1}_{\ell'}}.
\end{gather}

\begin{proposition}\label{asso}
The $\star$-product is associative, namely, for three arbitrary $q$-symmetric series
$\Ags(\bar u)$, $\Bgs(\bar u)$ and
$\Cgs(\bar u)$
 of the form~\eqref{gs5}
\begin{gather}\label{as-eq}
\sk{\Ags(\bar u)\star\Bgs(\bar u)}\star \Cgs(\bar u)=
\Ags(\bar u)\star\sk{\Bgs(\bar u)\star \Cgs(\bar u)}.
\end{gather}
\end{proposition}

\begin{proof}
 We check an equality \r{as-eq} f\/irst in the simplest case of the generating
series depending on one generating parameter $u$. Equating the coef\/f\/icients at the $n$-th
power of this parameter we obtain from \r{as-eq} an equality
\begin{gather}
\tSym_{\, \bar t}\sk{\sum_{n\geq m\geq 0}\sum_{n\geq s\geq m}
\Ags(t_{s+1},\ldots,t_{n})\cdot \Bgs(t_{m+1},\ldots,t_{s})\cdot \Cgs(t_{1},\ldots,t_m)
}\nonumber\\
\qquad{}=\tSym_{\, \bar t}\sk{\sum_{n\geq s'\geq 0}\sum_{s'\geq m'\geq 0}
\Ags(t_{s'+1},\ldots,t_{n})\cdot \Bgs(t_{m'+1},\ldots,t_{s'})\cdot \Cgs(t_{1},\ldots,t_{m'})
},\label{cas1}
\end{gather}
where the property \r{exa3} of the $q$-symmetric generating series was used and $\bar t$ is a set
$\{t_1,\ldots,t_n\}$.
An equality \r{cas1} is an obvious identity if one replaces the ordering of the summations.
It is clear that in the general case the arguments remain the same and the appearing of the series~\r{Zserm} does not change these arguments.\end{proof}

\begin{proposition}\label{inv-gs}
For any generating series $\Ags(\bar u)$ there exist an unique $q$-symmetric  series $\Bgs(\bar u)$
such that
\begin{gather*}
\Bgs(\bar u)\star\Ags(\bar u)=\Ags(\bar u) \star\Bgs(\bar u)=1.
\end{gather*}
\end{proposition}
\begin{proof}
Since  $\Ags(\bar u)$ has the form of a Taylor series with the free term equal to $1$,
we can always reconstruct uniquely the inverse series solving recursively the equations
for the coef\/f\/icients of the series $\Bgs(\bar u)$.
By the construction the coef\/f\/icients of this series will be also $q$-symmetric.
\end{proof}

\section[Universal nested Bethe vectors for $U_q(\widehat{\mathfrak{gl}}_N)$]{Universal nested Bethe vectors for $\boldsymbol{\Uqgln}$}

Quantum af\/f\/ine algebras in the current realization \cite{D88} provide examples
of the $q$-symmetric ge\-ne\-rating series.
We will construct these ge\-ne\-rating series for the quantum af\/f\/ine algebra
$\Uqgln$ and  show that  $\star$-products of these generating series
provide  hierarchical relations for~NBA. We now recall
the current realization of the algebra $\Uqgln$.

The quantum af\/f\/ine algebra $\Uqgln$ is generated by the modes of the currents
\begin{gather}\label{cur}
E_i(z)=\sum_{n\in\ZZ}E_i[n]z^{-n},\qquad
F_i(z)=\sum_{n\in\ZZ}F_i[n]z^{-n},\qquad
k^\pm_j(z)=\sum_{n\geq 0} k_j[\pm n]z^{\mp n},
\end{gather}
where $i=1,\ldots,N-1$ and $j=1,\ldots,N$
subject to the commutation relations given in the Appendix~\ref{com-rel}.
The generating series $F_i(z)$, $E_i(z)$ and $k^\pm(z)$ are called
{\em total and Cartan currents} respectively.

We consider two types of Borel subalgebras of the
algebra $\Uqgln$.

Generators of the standard Borel subalgebras $U_q(\mathfrak{b}^\pm)\subset \Uqgln $
can be expressed in terms of the modes of the currents \r{cur}. To do this one has to introduce
the composed currents~$E_{a,b}(z)$ and $F_{b,a}(z)$ for $a<b-1$ and $1\leq a<b\leq N$
(see Appendix~\ref{com-rel} for the
definition of the currents~$F_{b,a}(z)$).
The Borel subalgebra $U_q(\mathfrak{b}^+)$
is  generated by the modes of the currents:
$E_i[m]$, $m>0$; $F_i[n]$,  $k^+_j[n]$, $n\geq0$ and $E_{a,b}[1]$, $a<b-1$.
Dual standard Borel subalgebra $U_q(\mathfrak{b}^-)$
is  generated by the modes of the currents:
$F_i[m]$, $m<0$; $E_i[n]$,  $k^+_j[n]$, $n\leq0$ and $F_{b,a}[-1]$, $a<b-1$.
Here $i=1,\ldots,N-1$ and $j=1,\ldots,N$. The reader can find description of the standard Borel
subalgebras in terms of the modes of the $\Uqgl{3}$ currents  in the paper \cite{KP}.
This decomposition of the algebra $\Uqgln$ is related to the standard realization of this algebra
in terms of pair of the
 dual $\LL$-operators, where generators of the standard Borel subalgebras serve as the
modes of the Gauss coordinates of the corresponding $\LL$-operators.

Another type of Borel subalgebras is related to the current realization of
$\Uqgln$ and was introduced in \cite{D88}.
The Borel subalgebra $U_F\subset \Uqgln$ is generated by the modes
$F_i[n]$, $k^+_j[m]$, $i=1,\ldots,N-1$, $j=1,\ldots,N$, $n\in\ZZ$
and $m\geq0$. The Borel subalgebra
$U_E\subset \Uqgln$ is generated by the modes
$E_i[n]$, $k^-_j[-m]$, $i=1,\ldots,N-1$, $j=1,\ldots,N$, $n\in\ZZ$ and
$m\geq0$. We  also consider  a subalgebra $U'_F\subset U_F$, generated
by the elements
$F_i[n]$, $k^+_j[m]$, $i=1,\ldots,N-1$, $j=1,\ldots,N$, $n\in\ZZ$
and $m>0$, and a subalgebra $U'_E\subset U_E$ generated by
the elements
$E_i[n]$, $k^-_j[-m]$, $i=1,\ldots,N-1$, $j=1,\ldots,N$, $n\in\ZZ$
and $m>0$. We call these subalgebras of $\Uqgln$ the {\it current Borel subalgebras}.
Further, we will be interested in the intersections,
\begin{gather}
\label{Intergl}
U_f^-=U'_F\cap U_q(\mathfrak{b}^-),\qquad
U_F^+=U_F\cap U_q(\mathfrak{b}^+)
\end{gather}
and will describe properties of projections to these intersections.

The current Borel subalgebras are  Hopf subalgebras of $\Uqgln$ with respect
 to the \emph{current} Hopf structure for the algebra
$\Uqgln$ def\/ined in~\cite{D88}:
\begin{gather}
 \Delta^{(D)}\sk{E_i(z)} =E_i(z)\ot 1 + k^-_{i}(z)\sk{k^-_{i+1}(z)}^{-1}\ot
E_i(z) ,\nonumber\\
\Delta^{(D)}\sk{F_i(z)} =1\ot F_i(z) + F_i(z)\ot
k^+_{i}(z)\sk{k^+_{i+1}(z)}^{-1},\label{gln-copr}\\
\Delta^{(D)}\sk{k^\pm_i(z)} =k^\pm_i(z)\ot k^\pm_{i}(z) .\nonumber
\end{gather}
The quantum af\/f\/ine algebra $\Uqgln$ with ommited central charge and gradation operator
can be identif\/ied  with the quantum double of its current Borel subalgebra
constructed using the comultiplication \r{gln-copr}.

One may check that the intersections  $U_f^-$
and $U_F^+$ are subalgebras. It was proved in \cite{KPT}
that these subalgebras  are coideals with respect to Drinfeld coproduct
\r{gln-copr}
\begin{gather*}
\Delta^{(D)}(U_F^+)\subset \Uqgln\ot U_F^+,\qquad
\Delta^{(D)}(U_f^-)\subset U_f^-\ot \Uqgln,
\end{gather*}
and the multiplication $m$ in $\Uqgln$ induces an isomorphism of vector spaces
\begin{gather*}
m: \ \ U_f^-\ot U_F^+\to U_F.
\end{gather*}
According to the general theory presented in \cite{EKhP} we  def\/ine
 projection operators $\Pfp:U_F\subset \Uqgln \to U_F^+$ and
$\Pfm:U_F\subset \Uqgln \to U_f^-$ by the prescriptions
\begin{gather}\label{pgln}
\Pfp(\Ff_-\ \Ff_+)=\coun(\Ff_-)\ \Ff_+, \qquad
\Pfm(\Ff_-\ \Ff_+)=\Ff_-\ \coun(\Ff_+),
\qquad \text{for any}\  \Ff_-\in U_f^-,
\ \Ff_+\in U_F^+ ,\!\!\!\!
\end{gather}
where $\coun$ is the counit map: $\coun:\Uqgln\to\mathbb{C}$.

Denote by  $\overline U_F$ an extension of the algebra $U_F$ formed
by linear combinations of series, given as inf\/inite sums of monomials
$a_{i_1}[n_1]\cdots a_{i_k}[n_k]$ with $n_1\leq\cdots\leq n_k$, and $n_1+\dots +n_k$
f\/ixed,
where  $a_{i_l}[n_l]$ is either $F_{i_l}[n_l]$ or $k^+_{i_l}[n_l]$.
It was proved in \cite{EKhP} that
\begin{itemize}\itemsep=0pt
\item[(1)] the projections \r{pgln} can be extended to the
 algebra
$\overline U_F$;
\item[(2)] for any $\Ff\in \overline U_F$ with $\Delta^{(D)}(\Ff)=\sum_i \Ff'_i\otimes \Ff''_i$ we have
\begin{gather}\label{pr-prop}
\Ff=\sum_i\Pfm(\Ff''_i)\cdot \Pfp(\Ff'_i).
\end{gather}
\end{itemize}

\subsection{Generating series for universal Bethe vectors}

It was proved in the papers \cite{KPT,KhP-GLN} that the
projection of the product of the $\Uqgln$ currents can be identif\/ied with {\it universal
Bethe vectors} (UBV).
In this paper we show that the hierarchical relations for UBV
 can be presented in a compact form using  $\star$-product of certain
$q$-symmetric generating series.
Then the formal inversion of generating series allows to obtain another form of hierarchical
 relations and to investigate further (see~\cite{FKPR})  special
properties of UBV  when their parameters satisfy
the universal Bethe equations appeared in the framework of the analytical Bethe ansatz~\cite{ACDFR}. In this paper we will demonstrate these properties for the $\Uqgl{2}$ universal
Bethe vectors.

 Products of the $\Uqgln$ currents yield examples of the $q$-symmetric
generating series.
 We consider the generating series of the product of the currents
\begin{gather}\label{gs1}
\Fgs^N(\bar u)=\sum_{n_1,\ldots,n_{N-1}\geq0} \Fgs^N(\bar t_{\segg{\bar n}})
u_1^{n_1} u_2^{n_2}\cdots u_{N-1}^{n_{N-1}},
\end{gather}
where each term
$\Fgs^N(\bar t_{\segg{\bar n}})$ means the following normalized
product of the currents
\begin{gather}\label{gs2}
\Fgs^N(\bar t_{\segg{\bar n}})=\frac{F_{N-1}(t_{n_{N-1}}^{N-1})\cdots F_{N-1}(t_{1}^{N-1})
\cdots  F_{1}(t_{n_1}^1)\cdots F_{1}(t_{1}^1)}{n_{N-1}!\cdots n_{1}!}.
\end{gather}
We set $\Fgs^N_{\bar 0}\equiv 1$. More generally, following the convention \rf{anm},
for a segment $\meg{\bar\lll}{\bar\rr}$ and the
 related collection
$\bar t_{\meg{\bar\lll}{\bar\rr}}$  of variables, see \rf{t-coll}, we set
\begin{gather*}
\Fgs^N(\bar t_{\meg{\bar\lll}{\bar\rr}})=\frac{F_{N-1}(t_{r_{N-1}}^{N-1})\cdots F_{N-1}(t_{l_{N-1}+1}^{N-1})
\cdots   F_{1}(t_{r_1}^1)\cdots F_{1}(t_{l_1+1}^1)}{(r_{N-1}-l_{N-1})!\cdots (r_{1}-l_1)!}.
\end{gather*}

The superscript $N$ in the notation $\Fgs^N(\bar t_{\segg{\bar n}})$ of the coef\/f\/icients of the $q$-symmetric generating series
signifies that these coef\/f\/icients belong to the subalgebra $\U_F\subset\Uqgln$.
Further on we will consider smaller algebras
$\Uqgl{j+1}$, $j=1,\ldots,N-1$ embedded into $\Uqgln$ in two dif\/ferent way. Embedding
$\tau_j:\Uqgl{j+1}\hookrightarrow\Uqgln$ is def\/ined by removing the currents $F_{j+1}(t),\ldots,F_{N-1}(t)$,
$E_{j+1}(t),\ldots,E_{N-1}(t)$ and the Cartan currents $k^\pm_{j+2}(t),\ldots,k^\pm_{N}(t)$.
Embedding
$\tilde\tau_{j}:\Uqgl{j+1}\hookrightarrow\Uqgln$ is def\/ined by removing the currents $F_{1}(t),\ldots,F_{N-j-1}(t)$,
$E_{1}(t),\ldots,E_{N-j-1}(t)$ and the Cartan currents $k^\pm_{1}(t),\ldots,k^\pm_{N-j-1}(t)$.
Note, that the embedding $\tau_{N-1}=\tilde\tau_{N-1}$ is an identical map.

Using these embeddings, for any segment $\meg{\bar\lll}{\bar\rr}$ such that $l_i=r_i$ for all
$i=N-j+1,N-j+2,\ldots, N-1$ we def\/ine a series
\begin{gather*}
\Fgs^{N-j+1}(\bar t_{\meg{\bar\lll}{\bar\rr}})=\frac{F_{N-j}(t_{r_{N-j}}^{N-j})\cdots
F_{N-j}(t_{l_{N-j}+1}^{N-1})
\cdots   F_{1}(t_{r_1}^1)\cdots F_{1}(t_{l_1+1}^1)}{(r_{N-j}-l_{N-j})!\cdots (r_{1}-l_1)!}
\end{gather*}
and for any segment $\meg{\bar\lll}{\bar\rr}$ such that $l_i=r_i$ for all
$i=1,2,\ldots, j-1$  a series
\begin{gather*}
\tilde{\Fgs}^{N-j+1}(\bar t_{\meg{\bar\lll}{\bar\rr}})=\frac{F_{N-1}(t_{r_{N-1}}^{N-1})\cdots
F_{N-1}(t_{l_{N-1}+1}^{N-1})
\cdots   F_{j}(t_{r_j}^1)\cdots F_{j}(t_{l_j+1}^1)}{(r_{N-1}-l_{N-1})!\cdots (r_{j}-l_j)!} .
\end{gather*}
These series are gathered into $q$-symmetric generating functions $\Fgs^{N-j+1}(u_1,\ldots,u_{N-j})$ and
$\tilde{\Fgs}^{N-j+1}(u_j,\ldots,u_{N-1})$.
Subscripts $a=1,\ldots,N-1$ of the formal parameters $u_a$ in the de\-f\/i\-nition of these
generating series  denote the indices of the simple roots of the algebra
$\mathfrak{gl}_N$. For example, the notation $\Fgs^{N-1}(u_1,\ldots,u_{N-2})$ means the
generating series which coef\/f\/icients take values in the subalgebra $\Uqgl{N-1}$ embedded into
$\Uqgln$ by means of the map $\tau_{N-2}$. On the other hand,
 the notation $\tilde{\Fgs}^{N-1}(u_2,\ldots,u_{N-1})$ means the
generating series taking value in the subalgebra $\Uqgl{N-1}$ embedded into
$\Uqgln$ by means of the map $\tilde\tau_{N-2}$.

Further on we  use the special notation for the ordered product of the noncommutative entries.
Symbols $\mathop{\prod}\limits^{\longleftarrow}_a A_a$ and
$\mathop{\prod}\limits^{\longrightarrow}_a A_a$
 will mean  ordered products of
noncommutative entries~$A_a$, such that~$A_a$ is on the right (resp., on the left)
 from~$A_b$ for $b>a$:
\begin{gather*}
\mathop{\prod}\limits^{\longleftarrow}_{j\geq a\geq i} A_a = A_j A_{j-1}
\cdots  A_{i+1} A_i ,\qquad
\mathop{\prod}\limits^{\longrightarrow}_{i\leq a\leq j} A_a = A_i A_{i+1}
\cdots  A_{j-1} A_j .
\end{gather*}
Using these notations we write the $\bar n$-th term  of the generating series \r{gs1}
as follows
\begin{gather}\label{tFFFl}
\Fgs^N(\bar t_{\segg{\bar n}})=\prod_{N-1\geq a\geq
1}^{\longleftarrow} \frac{1}{n_a!} \sk{\prod_{n_a\geq\ell\geq 1}^{\longleftarrow} F_{a}(t^a_{\ell})} .
\end{gather}

Besides of the generating series of  products of the total currents we  consider
also the generating series of the projections of the products of the currents
\begin{gather}\label{gs3}
\Pfpm\sk{\Fgs^N(\bar u)}=\sum_{n_1,\ldots,n_{N-1}\geq0}
\Pfpm\sk{\Fgs^N(\bar t_{\segg{\bar n}})}\  u_1^{n_1}\,u_2^{n_2}\cdots u_{N-1}^{n_{N-1}},
\end{gather}
where
\begin{gather}\label{gs4}
\Pfpm\sk{\Fgs^N(\bar t_{\segg{\bar n}})}=
\frac{\Pfpm\sk{F_{N-1}(t_{n_{N-1}}^{N-1})\cdots F_{N-1}(t_{1}^{N-1})
\cdots F_{1}(t_{n_1}^1)\cdots F_{1}(t_{1}^1)}}{n_{N-1}!\cdots n_{1}!}.
\end{gather}
In the same manner we def\/ine series
\begin{gather*}
\Pfpm\sk{\Fgs^{N-j+1}(\bar u)}=
\Pfpm\sk{\Fgs^{N-j+1}(u_1,\ldots,u_{N-j})}
\end{gather*}
 and
\begin{gather*}
\Pfpm\sk{\tilde{\Fgs}^{N-j+1}(\bar u)}=\Pfpm\sk{\tilde{\Fgs}^{N-j+1}(u_j,\ldots,u_{N-1})}.
\end{gather*}

The series $\Pfp\sk{\Fgs^N(\bar u)}$ is the generating series of all possible
universal of\/f-shell Bethe vectors. Our goal is to show
that the hierarchical relations of the nested Bethe vectors imply the facto\-ri\-za\-tion
property of this generating series with respect to the $\star$-product of certain
$q$-symmetric generating series.
An associativity of this product allows to obtain a new presentation for the universal
Bethe vectors.

We call any expression $\sum_i f^{(i)}_-\cdot f^{(i)}_+$, where
$f^{(i)}_-\in U^-_f$ and $f^{(i)}_+\in U^+_F$ {\it (normal) ordered}.
\begin{proposition}\label{nogs}
The $q$-symmetric generating series \eqref{gs1}  can be written using a $\star$-product
in a~normal ordered form
\begin{gather}\label{n-o-gs}
\Fgs^N(\bar u)=\Pfm\sk{\Fgs^N(\bar u)}\star\Pfp\sk{\Fgs^N(\bar u)},
\end{gather}
where the $q$-symmetric
generating series $\Pfpm\sk{\Fgs^N(\bar u)}$ are defined by \eqref{gs3}.
\end{proposition}

\begin{proof} Using the property of the projections \r{pr-prop}
an equality of the series
\begin{gather}\label{dec-ff1m}
 \F^N(\bar t_{\segg{\bar n}})=
\sum_{0\leq \ss_{N-1}\leq n_{N-1}}\!\!\cdots\!\! \sum_{0\leq \ss_1\leq n_1}
 \tSym_{\ \bar t_{\segg{\bar n}}}
\left(Z_{\bar\ss}({\bar t}_{\segg{\bar n}})
  \Pfm\sk{\F^N(\bar t_{\seg{\bar n}{\bar\ss}})}\cdot
\Pfp\sk{\F^N(\bar t_{\seg{\bar\ss}{\bar0}})}\right)\!\!\!
\end{gather}
was proved in \cite{KhP-GLN} (see, Proposition~4.1 therein). That proof was based on the
 comultiplication  property \r{gln-copr} and the commutation relation between currents.
Using a def\/inition of the $\star$-pro\-duct and considering the coef\/f\/icients in front of monomial
$u_1^{n_1}\cdots u_{N-1}^{n_{N-1}}$ in the both sides of the equality~\r{n-o-gs} we
obtain the formal series equality \r{dec-ff1m}. \end{proof}

\section{Generating series of strings and nested Bethe ansatz}

\subsection{Composed currents and the strings}

For two sets of variables $\{t^1_1,\ldots,t^1_k\}$ and $\{t^2_1,\dots,t^2_k\}$
we introduce the series
\begin{gather}
  V(t^2_k,\ldots,t^2_1;t^1_k,\ldots,t^1_1) =  \prod_{m=1}^k\frac{t^1_m/t^2_m}{1-t^1_m/t^2_m}
\prod_{m'=m+1}^{k}\frac{q-q^{-1}t^1_{m'}/t^2_m}{1-t^1_{m'}/t^2_m}\nonumber\\
\phantom{V(t^2_k,\ldots,t^2_1;t^1_k,\ldots,t^1_1) }{} =  \prod_{m=1}^k\frac{t^1_m/t^2_m}{1-t^1_m/t^2_m}
\prod_{m'=1}^{m-1}\frac{q-q^{-1}t^1_m/t^2_{m'}}{1-t^1_m/t^2_{m'}}  .\label{rat-Y}
\end{gather}

Fix $j=1,\ldots,N-1$ and a collection of  non-negative integers
${\bar \ss}^j=\{\ss_1,\ldots,\ss_j\}$ satisfying the admissibility condition:
\begin{gather}\label{non-inc}
0=\ss_0\leq\ss_1\leq\ss_2\leq\cdots\leq \ss_{j-1}\leq\ss_{j} .
\end{gather}
 We def\/ine a series depending on the set of the variables
\begin{gather*}
 \bar t_{\segg{\bar\ss^j}}=\{t^1_{1},\ldots,t^1_{\ss_1}; t^2_{1},\ldots,t^2_{\ss_2}; \ldots;
t^j_{1},\ldots,t^j_{\ss_j}\}
 \end{gather*}
 of the form
\begin{gather}\label{rat-X}
X(\bar t_{\segg{\bar\ss^j}})=
\prod_{a=1}^{j-1}
V(t^{a+1}_{\ss_{a}},\ldots,t^{a+1}_{1};
t^{a}_{\ss_{a}},\ldots,t^{a}_{1} ) .
\end{gather}
When $j=1$ we set $X(\cdot)=1$.

Def\/ine an ordered normalized product of the composed currents, which we call
{\it a string} of the type $j$:\footnote{Here the def\/inition of the string dif\/fers from those used in
\cite{KhP-GLN,OPS} by the  combinatorial factor being a rational function of parameter $t_i^a$.}
\begin{gather}\label{string}
\Sgs^{j+1}(\bar t_{\segg{\bar\ss^j}})=X(\bar t_{\segg{\bar\ss^j}})
\prod^{\longleftarrow}_{j\geq a\geq 1}\sk{\frac{1}{ (\ss_a-\ss_{a-1})! }
\prod^{\longleftarrow}_{\ss_{a}\geq \ell>\ss_{a-1}}
F_{j+1,a}(t^j_\ell)}
\end{gather}
taking values in the subalgebra $\U_F\subset\Uqgl{j+1}$ embedded into $\Uqgln$ by the map $\tau_j$.
The composed currents $F_{j,i}(t)$   corresponding to  non-simple roots of the algebra
$\mathfrak{gl}_N$ belongs to the completion $\overline U_F$ and their def\/inition is given in the
Appendix~\ref{com-rel} (see \r{rec-f1}, \r{rec-f111}, \r{res-in}).

More generally,  let $\bar{m}=\{{m}_1,\ldots,{m}_{N-1}\}$ and
$\bar{n}=\{{n}_1,\ldots,{n}_{N-1}\}$ be a pair of collections of nonnegative
integers such that ${n}_a-{m}_a=0$ for $a=j+1,\ldots,N-1$ and $n_a-m_a=s_a$ for
any $a=1,\ldots, j$. Then for the set of variables
\begin{gather*}
 \bar t_{\segg{\bar{m},\bar{n}}}=\{t^1_{m_1+1},\ldots,t^1_{{n}_1}; t^{2}_{m_{2}},
 \ldots,t^{2}_{{n}_{2}};
 \ldots;  t^{j}_{m_{j}+1},\ldots,t^{j}_{{n}_{j}}\},
\end{gather*}
we set
\begin{gather*}
X(\bar t_{\segg{\bar{m},\bar{n}}})= \prod_{a=1}^{j-1}
V(t^{a+1}_{m_{a+1}+s_a},\ldots,t^{a+1}_{m_{a+1}+1};
t^{a}_{{n}_{a}},\ldots,t^{a}_{{m}_a+1} ),
\\
\Sgs^{j+1}
(\bar t_{\segg{\bar{m},\bar{n}}})=  X(\bar
t_{\segg{\bar{m},\bar{n}}}) \prod^{\longleftarrow}_{j\geq a \geq
1}\sk{\frac{1}{ ({s}_{a}-{s}_{a-1})! }
\prod^{\longleftarrow}_{s_a+m_j\geq\ell> s_{a-1}+m_j}
F_{j+1,a}(t^j_\ell)} .
\end{gather*}

We def\/ine the $q$-symmetric generating series of the strings of the type $j$ by the formula
\begin{gather}\label{gs-str}
\Sgs^{j+1}(u_1,\ldots,u_j)=\sum_{s_j\geq s_{j-1}\geq \cdots \geq s_1\geq 0}
\tSym_{\, \bar t_{\segg{\bar\ss^j}}}\sk{
\Sgs^{j+1}(\bar t_{\segg{\bar\ss^j}})}  u_1^{s_1} u_2^{s_2} \cdots
u_{j-1}^{s_{j-1}}u_j^{s_j} .
\end{gather}
Here superscript $j+1$ signif\/ies that this generating series takes values in
the subalgebra $\Uqgl{j+1}$ embedded into $\Uqgln$ by the map $\tau_j$.
The subscripts of the  parameters $u_1,\ldots,u_j$ signify that this
 subalgebra is generated by the $\Uqgln$ currents corresponding
to the simple roots with indices $1,\ldots,j$.

Let $\Pfp\sk{\Fgs^N(u_1,\ldots,u_{N-1})}$ be the generating series of the universal of\/f-shell Bethe vectors
for the algebra $\Uqgln$ and $\Pfp\sk{\Fgs^{N-1}(u_1,\ldots,u_{N-2})}$ be the analogous series for the
smaller algebra $\Uqgl{N-1}$ embedded into $\Uqgln$ by the map $\tau_{N-2}$.  We have the following
\begin{proposition}\label{NBA}
Hierarchical relations between universal weight functions for algebras $\Uqgln$ and
$\Uqgl{N-1}$ can be written as the following equality on
generating series:
\begin{gather}\label{rrNBA}
\Pfp\sk{\Fgs^N(u_1,\ldots,u_{N-1})}=\Pfp\sk{\Sgs^{N}(u_1,\ldots,u_{N-1})} \star
\Pfp\sk{\Fgs^{N-1}(u_1,\ldots,u_{N-2})} .
\end{gather}
\end{proposition}

\begin{proof}Taking the coef\/f\/icients in front of the monomial
$u_1^{n_1}u_2^{n_2}\cdots u_{N-1}^{n_{N-1}}$ we obtain an equality of the formal series
\begin{gather}
\Pfp\sk{\F^{N}(\bar t_{\segg{\bar n}})} =\sum_{n_{N-1}=s_{N-1}\geq \cdots \geq s_1\geq 0}
\tSym_{\, \bar t_{\segg{\bar n}}} \Big(
Z_{\bar s}(\bar t_{\segg{\bar n}})\nonumber\\
\phantom{\Pfp\sk{\F^{N}(\bar t_{\segg{\bar n}})} =}{}
\Pfp\sk{\Sgs^{N}(\bar t_{\seg{\bar n}{\bar n-\bar\ss}})}
\cdot \Pfp\sk{\F^{N-1}(\bar t_{\segg{{\bar n}-\bar\ss}})}\Big) .\label{rrNBA1}
\end{gather}
An equality \r{rrNBA1}
coincides with the statement of the Proposition~4.2 of the paper \cite{KhP-GLN} up to renormalization
of the universal weight function by the combinatorial factor.
\end{proof}

\begin{corollary}
Generating series of the $\Uqgln$ universal weight functions can be written using ordering
$\star$-product of the generating series of the strings
\begin{gather}
\Pfp\sk{\Fgs^N(u_1,\ldots,u_{N-1})} =\Pfp\sk{\Sgs^{N}(u_1,\ldots,u_{N-1})} \star
\Pfp\sk{\Sgs^{N-1}(u_1,\ldots,u_{N-2})} \star
\cdots\nonumber\\
 \phantom{\Pfp\sk{\Fgs^N(u_1,\ldots,u_{N-1})} =}{} \star  \Pfp\sk{\Sgs^{3}(u_1,u_{2})} \star  \Pfp\sk{\Sgs^{2}(u_{1})}  .\label{rrNBA2}
\end{gather}
\end{corollary}

 In \r{rrNBA2} we assume that the universal Bethe vectors for the algebra $\Uqgl{1}$ are equal~to~1.

Recall that the generating series $\Pfp\sk{\Sgs^{j}(u_1,\ldots,u_{j-1})}$ belongs to subalgebra
$\overline U_F\subset \Uqgl{j}$ embedded into $\Uqgln$ by the map $\tau_{j-1}$ which removes the currents
$F_a(t)$, $E_a(t)$ and $k^\pm_{a+1}(t)$ with $a=j,j+1,\ldots,N-1$.

\subsection{Hierarchical relations for the negative projections}

One of the results of the papers \cite{KhP-GLN,OPS} is the hierarchical relation
for the positive projections of the product of the currents. Let us recall shortly the main idea of
this calculation. In order to calculate the projection
\begin{gather*}
\Pfp\sk{
F_{N-1}(t^{N-1}_{n_{N-1}})\cdots F_{N-1}(t^{N-1}_{1})
\cdots
   F_{1}(t^{1}_{n_{1}})\cdots F_{1}(t^{1}_{1})   }
\end{gather*}
we separate all factors $F_a(t^a_\ell)$ with $a<N-1$ and apply to this product the
ordering procedure based on the property \r{dec-ff1m}. We obtain under total projection
the $q$-symmetrization of terms $x_i\Pfm(y_i)\Pfp(z_i)$, where $x_i$ are expressed
via modes of $F_{N-1}(t)$ and $y_i$, $z_i$ via modes of $F_a(t)$ with $a<N-1$.
Then we used the property of the projection that
\begin{gather*}
\Pfp\sk{x_i\Pfm(y_i)\Pfp(z_i)}=\Pfp\sk{x_i\Pfm(y_i)}\cdot \Pfp(z_i)
\end{gather*}
and reorder the product of $x_i$ and $\Pfm\sk{y_i}$ under positive projection to obtain
the string build from the composed currents (cf. equation~\r{rrNBA1}).

We will use an analogous strategy to calculate the negative projection of the same product of the
currents
\begin{gather}\label{nuWF}
\Pfm\sk{
F_{N-1}(t^{N-1}_{n_{N-1}})\cdots F_{N-1}(t^{N-1}_{1})  \cdots
F_{2}(t^{2}_{n_{2}})\cdots
F_2(t^2_{1})  F_{1}(t^{1}_{n_{1}})\cdots F_{1}(t^{1}_{1})   } .
\end{gather}
Now we separate all factors $F_a(t^a_\ell)$ with $a>1$ and apply to this product the
ordering rule~\r{dec-ff1m}. Again, we obtain under total negative projection
the $q$-symmetrization of terms $\Pfm(x_i)\Pfp(y_i)z_i$, where $z_i$ are expressed
via modes of $F_{1}(t)$ and $y_i$, $z_i$ via modes of $F_a(t)$ with $a>1$.
Using now the property of the projections
\begin{gather*}
\Pfm\sk{\Pfm(x_i)\Pfp(y_i)z_i}=\Pfm\sk{x_i}\cdot \Pfm\sk{\Pfp(y_i) z_i}
\end{gather*}
and reordering the product of $\Pfp(y_i)$ and  $z_i$ under negative projection we
obtain the desired hierarchical relations for the negative projection of the currents product.

We will not repeat these calculations since they are analogous to the ones presented in
\cite{KhP-GLN}, but will formulate the f\/inal answer of these hierarchical relations.
For two sets of variables $\{t^1_1,\ldots,t^1_l\}$ and $\{t^2_1,\dots,t^2_l\}$
we introduce the series
\begin{gather}
  \tilde V(t^2_l,\ldots,t^2_1;t^1_l,\ldots,t^1_1) =\ds \prod_{m=1}^l\frac{1}{1-t^1_m/t^2_m}
\prod_{m'=m+1}^{l}\frac{q-q^{-1}t^1_{m'}/t^2_m}{1-t^1_{m'}/t^2_m}\nonumber\\
\phantom{\tilde V(t^2_l,\ldots,t^2_1;t^1_l,\ldots,t^1_1)}{}  =  \prod_{m=1}^l\frac{1}{1-t^1_m/t^2_m}
\prod_{m'=1}^{m-1}\frac{q-q^{-1}t^1_m/t^2_{m'}}{1-t^1_m/t^2_{m'}}  .\label{rat-Yt}
\end{gather}

Fix $k=1,\ldots,N-1$ and collection of the non-negative integers ${\tilde
{s}}^{{k}}=\{{s}_k,\ldots,{s}_{N-1}\}$ satisfying the admissibility condition
\begin{gather}\label{m-non-inc}
{s}_k\geq{s}_{k+1}\geq{s}_{k+2}\geq\cdots\geq {s}_{N-1}\geq{s}_{N}=0 .
\end{gather}
 We def\/ine a series depending on the set of the variables
\begin{gather*}
 \bar t_{\segg{\tilde{s}^{{k}}}}=\{t^k_{1},\ldots,t^k_{{s}_k}; t^{k+1}_{1},\ldots,t^{k+1}_{{s}_{k+1}};
 \ldots ; t^{N-1}_{1},\ldots,t^{N-1}_{{s}_{N-1}}\}
 \end{gather*}
 of the form
\begin{gather}\label{t-rat-X}
\tilde X(\bar t_{\segg{\tilde{s}^{{k}}}})= \prod_{a=k}^{N-2} \tilde
V(t^{a+1}_{{s}_{a+1}},\ldots,t^{a+1}_{1};
t^{a}_{{s}_{a}},\ldots,t^{a}_{{s}_a-{s}_{a+1}+1} ).
\end{gather}
When $k=N-1$ we set $\tilde X(\cdot)=1$.

Def\/ine an ordered normalized product of the composed currents, which we call
{\it a dual string} of the type $N-k$:
\begin{gather}\label{dstring}
\tSgs^{N-k+1}
(\bar t_{\segg{\tilde{s}^{{k}}}})=\tilde X(\bar
t_{\segg{\tilde{s}^{{k}}}}) \prod^{\longleftarrow}_{N\geq a > k}\sk{\frac{1}{
({s}_{a-1}-{s}_{a})! } \prod^{\longleftarrow}_{{s}_k-{s}_{a}\geq
\ell>{s}_k-{s}_{a-1}} F_{a,k}(t^k_\ell)} .
\end{gather}
Note that the notion of the dual string is dif\/ferent from the notion of the inverse string
used in~\cite{KhP-GLN}.

More generally,  let $\bar{m}=\{{m}_1,\ldots,{m}_{N-1}\}$ and
$\bar{n}=\{{n}_1,\ldots,{n}_{N-1}\}$ be a pair of collections of nonnegative
integers such that $n_a-m_a=0$ for $a=1,\ldots,k-1$ and
$n_a-m_a=s_a$ for
any $a=k,\ldots, N-1$. Then for the collection of variables
\begin{gather*}
 \bar t_{\segg{\bar{m},\bar{n}}}=\big\{t^k_{m_k+1},\ldots,t^k_{{n}_k};
  t^{k+1}_{m_{k+1}},
 \ldots,t^{k+1}_{{n}_{k+1}};
 \ldots;  t^{N-1}_{m_{N-1}+1},\ldots,t^{N-1}_{{n}_{N-1}}\big\},
\end{gather*}
we set
\begin{gather*}
\tilde X(\bar t_{\segg{\bar{m},\bar{n}}})=\prod_{a=k}^{N-2} \tilde
V(t^{a+1}_{{n}_{a+1}},\ldots,t^{a+1}_{m_{a+1}+1};
t^{a}_{{n}_{a}},\ldots,t^{a}_{{n}_a-{s}_{a+1}+1} ),
\\
\tSgs^{N-k+1}
(\bar t_{\segg{\bar{m},\bar{n}}})= \tilde X(\bar
t_{\segg{\bar{m},\bar{n}}}) \prod^{\longleftarrow}_{N\geq a >
k}\sk{\frac{1}{ ({s}_{a-1}-{s}_{a})! }
\prod^{\longleftarrow}_{{n}_k-{s}_{a}\geq \ell>{n}_k-{s}_{a-1}}
F_{a,k}(t^k_\ell)} .
\end{gather*}

Doing the calculations described above we obtain the recurrence relations for the negative
projections
\begin{gather}\label{mrr}
\Pfm\sk{\Fgs^{N}(\bar t_{\segg{\bar n}})}=\!\!\sum_{n_{1}=\!\ss_{1}\geq \cdots \geq \ss_{N-1}\geq 0}
\!\!\!\tSym_{\, \bar t_{\segg{\bar n}}}\!
\sk{Z_{\bar\ss}(\bar t_{\segg{\bar n}})\cdot\! \Pfm\!\sk{\Fgs^{N-1} (\bar t_{\seg{\bar n}{\bar\ss}})}
\cdot \!\Pfm\!\sk{\tSgs^{N}(\bar t_{\segg{\bar\ss}})}
}\!.\!\!\!\!
\end{gather}

 We def\/ine the generating series of the dual strings of the type $N-j$ by the formula
\begin{gather*}
\tSgs^{N-k+1}(u_k,\ldots,u_{N-1})=\!\!\sum_{\ss_k\geq \ss_{k+1}\geq \cdots \geq \ss_{N-1}\geq 0}
\!\!\tSym_{\,\bar t_{\segg{\tilde\ss^{k}}}}\sk{
\tSgs^{N-k+1}(\bar t_{\segg{\tilde\ss^{k}}})}  u_k^{\ss_k} u_{k+1}^{\ss_{k+1}} \cdots
u_{N-1}^{\ss_{N-1}}
\end{gather*}
taking values in the subalgebra $\overline U_F\subset \Uqgl{N-k+1}$ embedded into $\Uqgln$ by the
map $\tilde\tau_{N-k}$ which removes the  currents $F_a(t)$, $E_a(t)$ and $k^\pm_{a}(t)$
for $a=1,\ldots,k-1$.
The recurrence relations~\r{mrr} can be written as the $\star$-product of the generating
series
\begin{gather}\label{mrr1}
\Pfm\sk{\Fgs^N(u_1,\ldots,u_{N-1})}=
\Pfm\sk{\tilde\Fgs^{N-1}(u_2,\ldots,u_{N-1})}  \star    \Pfm\sk{\tSgs^{N}(u_1,\ldots,u_{N-1})} .
\end{gather}

Generating series of the negative projections of the product of the currents can be written using ordering
$\star$-product of the generating series of the dual strings
\begin{gather}
\Pfm\sk{\Fgs^N(u_1,\ldots,u_{N-1})} =
\Pfm\sk{\tSgs^{2}(u_{N-1})} \star
\Pfm\sk{\tSgs^{3}(u_{N-2},u_{N-1})} \star \cdots\nonumber\\
\phantom{\Pfm\sk{\Fgs^N(u_1,\ldots,u_{N-1})} =}{}  \star \Pfm\sk{\tSgs^{N-1}(u_2,\ldots,u_{N-1})} \star
\Pfm\sk{\tSgs^{N}(u_1,\ldots,u_{N-1})}
 .\label{mrr2}
\end{gather}

\subsection{Other type of the hierarchical relations}

A special ordering property of the universal Bethe vectors when their parameters
$t^a_\ell$ satisfy the universal Bethe equations~\cite{ACDFR} was investigated in~\cite{FKPR}. This property
leads to the fact that the ordering of the product of the universal transfer matrix
and the universal nested Bethe vectors is proportional to the {\it same} Bethe vector modulo the terms
which belong to some ideal in the algebra
if the parameters of this vector satisfy the universal Bethe equations. We will demonstrate
this property for the $\Uqgl{2}$ universal Bethe vectors in the Section~\ref{gl2-ex}.

A cornerstone of this ordering property lies in a new hierarchical relations for the universal
Bethe vectors, which can be proved using the technique of the generating series. Here we
give the detailed proof of the relation which particular form was used in the paper \cite{FKPR}.

Using normal ordering relation \r{n-o-gs} and \r{mrr1}
 we may write the generating series of the product of the currents in the form
 \begin{gather}
\Fgs^N(u_1,\ldots,u_{N-1})=\Pfm\sk{\Fgs^N(u_1,\ldots,u_{N-1})} \star
\Pfp\sk{\Fgs^N(u_1,\ldots,u_{N-1})} \label{cal1}\\
\qquad{}=
\Pfm\sk{\tilde\Fgs^{N-1}(u_2,\ldots,u_{N-1})} \star    \Pfm\sk{\tSgs^{N}(u_1,\ldots,u_{N-1})}
 \star \Pfp\sk{\Fgs^N(u_1,\ldots,u_{N-1})} .\nonumber
\end{gather}
On the other hand these generating series may be presented as the factorized product
\begin{gather}\label{cal2}
\Fgs^N(u_1,\ldots,u_{N-1})=\tilde\Fgs^{N-1}(u_2,\ldots,u_{N-1})\cdot \Fgs^2(u_1)=
\tilde\Fgs^{N-1}(u_2,\ldots,u_{N-1}) \star  \Fgs^2(u_1) .
\end{gather}
Applying the ordering relation \r{n-o-gs} to the series $\tilde\Fgs^{N-1}(u_2,\ldots,u_{N-1})$  again
we obtain an alternative to \r{cal1} expression for the generating series
$\Fgs^N(u_1,\ldots,u_{N-1})$:
 \begin{gather}
 \Fgs^N(u_1,\ldots,u_{N-1})=\tilde\Fgs^{N-1}(u_2,\ldots,u_{N-1}) \star \Fgs^2(u_1)\nonumber\\
 \qquad=
\Pfm\sk{\tilde\Fgs^{N-1}(u_2,\ldots,u_{N-1})}  \star   \Pfp\sk{\tilde\Fgs^{N-1}(u_2,\ldots,u_{N-1})}
\star \Fgs^2(u_1) .\label{cal3}
\end{gather}
Equating the right hand sides of \r{cal1} and \r{cal3} we obtain the identity
\begin{gather}
\Pfm\!\sk{\tSgs^{N}(u_1,\dots,u_{N-1})}
\star\Pfp\sk{\Fgs^{N}(u_1,\dots,u_{N-1})}\nonumber\\
\qquad{}=
\Pfp\!\sk{\tilde\Fgs^{N-1}(u_2,\dots,u_{N-1})}
\star\Fgs^2(u_1)\label{usid1}
\end{gather}
or
\begin{gather}
\Pfp \sk{\Fgs^N(u_1,\ldots,u_{N-1})}\nonumber\\
\qquad{}=\sk{\Pfm\sk{\tSgs^{N}(u_1,\ldots,u_{N-1})}}^{-1}
 \star
\Pfp\sk{\Fgs^{N-1}(u_2,\ldots,u_{N-1})}
\star \Fgs^2(u_1) .\label{usid}
\end{gather}

The identity \r{usid} relates the universal of\/f-shell Bethe vectors for the algebra
$\Uqgln$ and for the smaller algebra $\Uqgl{N-1}$. They can be considered
as an universal formulation of the relation used in the pioneer paper \cite{KR83} for the
obtaining the nested Bethe equations.
The equality~\r{usid} between generating series contains many hierarchical relations between
UBV. In order to get some particular identities between these UBV one has to invert
explicitly the generating series $\Pfm\sk{\tSgs^{N}(u_1,\ldots,u_{N-1})}$. This will be done
in the next subsection.

\subsection{Inverting generating series of the strings}

Let $\Ccal(\bar u)$ be the generating series of negative projections of dual strings of the type
$N-1$:
$\Ccal(\bar u)=\Pfm\sk{\tSgs^{N}(u_1,\ldots,u_{N-1})}$ and $\Dcal(\bar u)$ be the inverse series:
$\Dcal(\bar u) \star \Ccal(\bar u)=1$.
By  the construction, see \rf{dstring}, the coef\/f\/icients
\begin{gather}\label{coef-m}
\Ccal(\bar t_{\segg{\bar\mm}})
=\Pfm\sk{\tSym_{\, \bar t_{\segg{\bar\mm}}}
\sk{\tSgs^N(\bar t_{\segg{\bar\mm}})}}
\end{gather}
 of the generating series $\Ccal(\bar u)$
 are nonzero only if the admissibility conditions $m_1\geq m_2\geq\cdots \geq m_{N-1}\geq 0$
for the set $\{\bar\mm\}=\{ m_1,m_{2},\ldots,m_{N-1}\}$ are satisf\/ied. We have to f\/ind
coef\/f\/icients of the generating series $\Dcal(\bar u)$
such that $\Dcal(\bar u) \star \Ccal(\bar u)=1$.
The latter equality is equivalent to the system of equations
\begin{gather}\label{rec1}
\tSym_{\, \bar t_{\segg{\bar n}}}
\Big(
\mathop{\sum_{n_1\geq m_1\geq 0}\cdots \sum_{n_{N-1}\geq m_{N-1}\geq 0}}
\limits_{m_1\geq \cdots \geq m_{N-1}}
 Z_{\bar m}(\bar t_{\segg{\bar n}})
\Dcal(\bar t_{\seg{\bar n}{\bar m}})
\cdot \Ccal(\bar t_{\segg{\bar m}})\Big) =0.
\end{gather}
for the unknown functional coef\/f\/icients $\Dcal(\bar t_{\segg{\bar k}})$
 for all possible f\/ixed values
of $n_1,\ldots,n_{N-1}$.
\begin{proposition}\label{inversion}
 The coefficients $\Dcal(\bar t_{\segg{\bar k}})$ are nonzero only if admissibility
conditions $k_1\geq k_2\geq\cdots \geq k_{N-1}\geq 0$ are satisfied.
\end{proposition}

\begin{proof} can be performed recursively by considering f\/irst the cases for
$n_1>0$ and $n_2=\cdots=n_{N-1}=0$.
  For these values
of $\bar n$ the series $Z_{\bar m}(\bar t_{\segg{\bar n}})=1$ and coef\/f\/icients of the
series in \r{rec1} depend only on the variables of the f\/irst
type $t_1^1,\ldots t_{n_1}^1$. The relation \r{rec1} takes in this case the form
\begin{gather}\label{rec2}
\tSym_{\, \bar t_{\segg{n}}}
\sk{\sum_{n_1\geq m_1\geq0}
\Dcal(t^1_{m_1+1},\ldots,t^1_{n_1})
\cdot \Ccal(t^1_{1},\ldots,t^1_{m_1})}=0 ,
\end{gather}
Its solution can be written in the form
\begin{gather}\label{rec3}
\Dcal(\bar t_{\segg{n}})
=
\tSym_{\, \bar t_{\segg{n}}}
\left(\sum_{p=0}^{n_1-1}(-1)^{p+1}
\!\!\!\!\sum_{n_1=k_{p+1}>k_{p}>\cdots>k_1>k_0=0}
\prod^{\longleftarrow}_{p+1\geq r\geq 1}
 \Ccal(t^1_{k_{r-1}+1},\ldots,t^1_{k_r})\right).
 \end{gather}
After this we consider the relation \r{rec1} for arbitrary $n_1\!\geq0$, $n_2\!=1$ and
\mbox{$n_3=\!\cdots\!= n_{N-1}\!=0$}. Avoiding writing the dependence on the `$t$' parameters, that is using
notations $\Ccal_{\bar k}$ instead of~$\Ccal(t_{\segg{\bar k}})$, the relation
\r{rec1} takes the form
\begin{gather}
\tSym_{\, \bar t_{\segg{\bar n}}}
\sk{\sum_{m_1=0}^{n_1}   \Dcal_{n_1-m_1,1,\mathbf{0}}
\cdot \Ccal_{m_1,0,\mathbf{0}}} \nonumber\\
\qquad{} +
\tSym_{\, \bar t_{\segg{\bar n}}}
\sk{\sum_{m_1=1}^{n_1}    Z_{m_1,1,\mathbf{0}}\cdot \Dcal_{n_1-m_1,0,\mathbf{0}}
\cdot \Ccal_{m_1,1,\mathbf{0}}  }
=0 .\label{rec4}
\end{gather}
The rational series $Z_{m_1,0,\mathbf{0}}$ disappears in the f\/irst sum of \r{rec4} by the same
reason as in \r{rec2}.
Let us consider the relation \r{rec4} for $n_1=0$. The second sum is absent and the f\/irst
sum contains
only one terms $\Dcal_{0,1,\mathbf{0}} \cdot \Ccal_{0,0,\mathbf{0}}=0$ which is equal to zero.
This proves that the coef\/f\/icient $\Dcal_{0,1,\mathbf{0}} =0$ vanishes identically.
Now the f\/irst sum in the relation \r{rec4} is terminated at $m_1=n_1-1$ and
this relation allows to f\/ind all coef\/f\/icients $\Dcal_{n_1,1,\mathbf{0}}$ starting from
$\Dcal_{1,1,\mathbf{0}}=-Z_{1,1,\mathbf{0}}\Ccal_{1,1,\mathbf{0}}$. Considering the
relation \r{rec1} for
$n_1\geq0$, $n_2=2$ and $n_3=\cdots=n_{N-1}=0$ we prove f\/irst that
$\Dcal_{0,2,\mathbf{0}}=\Dcal_{1,2,\mathbf{0}} =0$ and then can f\/ind all coef\/f\/icients
$\Dcal_{n_1,2,\mathbf{0}}$ starting from
$\Dcal_{2,2,\mathbf{0}}$.
It is clear now that the coef\/f\/icients $\Dcal_{n_1,n_2,\mathbf{0}}$ are non-zero only
if $n_1\geq n_2$. Continuing we prove the statement of the proposition.
\end{proof}

\subsection{Inversion and combinatorics}\label{combi}

To invert explicitly the generating series of the projection of the strings we have to introduce
certain combinatorial data. First of all, according to Proposition~\ref{inversion}
we f\/ix a sequence of non-negative integers $n_1\geq n_2\geq\cdots\geq n_{N-1}$ and the
corresponding set of the variables $\bar t_{\segg{\bar n}}$.

\begin{figure}
\begin{center}
\begin{picture}(100,60)
\unitlength=1mm
\put(-5,1){\footnotesize1}
\put(-5,6){\footnotesize2}
\put(-5,11){\footnotesize3}
\put(-5,16){\footnotesize4}
\put(0,15){\line(1,0){25}}
\put(0,10){\line(1,0){25}}
\put(0,5){\line(1,0){30}}
\put(0,20){\line(1,0){10}}
\put(0,0){\line(1,0){30}}
\put(0,0){\line(0,1){20}}
\put(5,0){\line(0,1){20}}
\put(10,0){\line(0,1){20}}
\put(15,0){\line(0,1){15}}
\put(20,10){\line(0,1){5}}
\put(25,10){\line(0,1){5}}
\put(20,0){\line(0,1){5}}
\put(25,0){\line(0,1){5}}
\put(30,0){\line(0,1){5}}
\end{picture}
\end{center}
\vspace{-5mm}
\caption{Example of the connected diagram for $n=16$, $p=4$, $\chi_1=6$,
$\chi_2=3$, $\chi_3=5$, $\chi_4=2$.}
\label{fig1}
\end{figure}
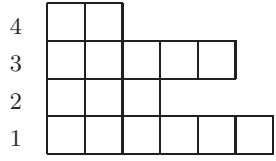

Choose any positive integer $n$ and
$p=1,\ldots,n$. A {\it diagram} $\chi$  of {\it size} $|\chi|=n$ and {\it height}
$p=\heit(\chi)$
is an
\underline{ordered} decomposition of $n$ into a sum of $p$ nonnegative integers,
\begin{gather*}
n=\chi_1+\cdots+\chi_{p}.
\end{gather*}
Equivalently,
 a diagram $\chi$  consists of $p$ rows and the $i$-th row contains $\chi_i$ boxes,
  An example of such diagram is
shown in the Fig.~\ref{fig1}. The rows of the diagrams are numbered from the bottom to
the top.
If $\chi_i=0$ for some $i=1,\ldots,p$ then  the diagram contains several
disconnected pieces. We will call the diagram $\chi$ {\it connected} if all
$\chi_i\not=0$ for $i=1,\ldots,p$.

A {\it tableaux} $\tchi$ with a given diagram $\chi$ is a f\/illing of all boxes of $\chi$ by
the indices $\{1,2,\dots,N-1\}$ of the positive roots of the algebra $\mathfrak{gl}_N$
with the condition of non-increasing from the left to the right along the rows.
If an index $a$ is associated to a box of the tableaux, we  say that this box has a `type' $a$.
We will call  the tableaux associated to the connected
diagrams {\it the  connected tableaux}.
For a given tableaux we def\/ine its weight $\wt(\tchi)$ as a set of numbers
$\bar n=\bar n(\tchi)=\{n_1(\tchi),\ldots,n_{N-1}(\tchi)\}$ such that $n_a(\tchi)$
is a number of boxes which have type bigger or equal than $a$, $a=1,\ldots,N-1$.
The size and  the height  of the tableaux is the size and
 the height of the corresponding diagram.
An example of a connected tableaux
is given on the Fig.~\ref{fig11}.

\begin{figure}
\begin{center}
\begin{picture}(100,60)
\unitlength=1mm
\put(0,15){\line(1,0){25}}
\put(0,10){\line(1,0){25}}
\put(0,5){\line(1,0){30}}
\put(0,20){\line(1,0){10}}
\put(0,0){\line(1,0){30}}
\put(0,0){\line(0,1){20}}
\put(5,0){\line(0,1){20}}
\put(10,0){\line(0,1){20}}
\put(15,0){\line(0,1){15}}
\put(20,10){\line(0,1){5}}
\put(25,10){\line(0,1){5}}
\put(20,0){\line(0,1){5}}
\put(25,0){\line(0,1){5}}
\put(30,0){\line(0,1){5}}
\put(2,1){\footnotesize$2$}
\put(7,1){\footnotesize$2$}
\put(12,1){\footnotesize$2$}
\put(17,1){\footnotesize$1$}
\put(22,1){\footnotesize$1$}
\put(27,1){\footnotesize$1$}
\put(2,16){\footnotesize$3$}
\put(7,16){\footnotesize$1$}
\put(2,6){\footnotesize$1$}
\put(7,6){\footnotesize$1$}
\put(12,6){\footnotesize$1$}
\put(2,11){\footnotesize$3$}
\put(7,11){\footnotesize$3$}
\put(12,11){\footnotesize$2$}
\put(17,11){\footnotesize$2$}
\put(22,11){\footnotesize$1$}
\end{picture}
\end{center}
\vspace{-5mm}

\caption{Example of connected tableaux in case of $N=4$ and its weight
$\{n_1,n_2,n_3\}=\{16,8,3\}$.}
\label{fig11}
\end{figure}

Denote by $\tchi^i$
the $i$th row of the tableaux $\tchi$.
Denote by $c^i_a=c^i_a(\tchi)$  the number of type $a$ boxes in the row $\tchi^i$.
Set $d^i_a=d^i_a(\tchi)=c^i_{N-1}+\cdots+c^i_{a}$ and
 $h^i_a=d^1_a+\cdots+d^i_a$
In particular, $d^i_1(\tchi)$ is the length $\chi_i$ of the row $\tchi^i$ and the collection
${\bar d}^i=\{d^i_1,\dots, d^i_{N-1}\}$ is the character
$\bar n(\tchi^i)$ of the row $\tchi^i$ considered as a tableaux by itself.
Clearly $d_a^i\geq d_b^i$ when $a\leq b$ and
\begin{gather*}
n_a=h_a^p=d_a^1+\cdots +d_a^p
\end{gather*}
for all $a=1,\ldots N-1$.
This formula demonstrates in particular that the weight of a tableaux always satisf\/ies
the admissibility conditions $n_1\geq n_2\geq \cdots\geq n_{N-1}$.


To each connected tableaux $\tchi$ of the  weight $\bar n=\bar n(\tchi)$
we associate a decomposition of the set of the variables $\bar t_{\segg{\bar n}}$ into the union
of $|\chi|$ disjoint subsets, each corresponding to a box of the diagram $\chi$ of tableaux $\tchi$.
To each box of the type $a$ we associate one variable of the type~1, one variables of the type~2, etc.,
one variable of the type~$a$, altogether $a$ variables.
We will number variables of the each type starting from the most bottom row and the most right box
where variable of this type appear for the f\/irst time.

Let us give an example of the decomposition and of the ordering for the tableaux
shown on the Fig.~\ref{fig3}.
The most bottom and the right box has the type 1. We associate to this box one
variable\footnote{Recall that superscript $a$ signif\/ies the `type' of the variable
$t^a_i$ and subscript $i$ counts the number of the variables of this type.}
$t^1_1$ of the
same type. Next to the left along the same row box has the type 2. We associate to this box two variables
$t^1_2$ and $t^2_1$ of the types 1 and 2. Last box in the bottom row has type 3 and we associate to this
box three variables $t^1_3$, $t^2_2$ and $t^3_1$. Next box is in the next row and also has the type $3$.
To this box we associate also three variables $t^1_4$, $t^2_3$ and $t^3_2$. Next two boxes in the third row
both have the type 1 and we associate to the most right box in this row one variable $t^1_5$ and to the
last box also one variable $t^1_6$.

\begin{figure}
\begin{center}
\begin{picture}(230,70)
\unitlength=1.5mm
\put(0,0){\line(0,1){15}}
\put(5,0){\line(0,1){15}}
\put(15,0){\line(0,1){5}}
\put(10,10){\line(0,1){5}}
\put(10,0){\line(0,1){5}}
\put(0,0){\line(1,0){15}}
\put(0,5){\line(1,0){15}}
\put(0,10){\line(1,0){10}}
\put(0,15){\line(1,0){10}}
\put(2,1.5){$3$}
\put(7,1.5){$2$}
\put(12,1.5){$1$}
\put(2,6.5){$3$}
\put(2,11.5){$1$}
\put(7,11.5){$1$}
\put(30,0){\line(0,1){15}}
\put(42,0){\line(0,1){15}}
\put(57,0){\line(0,1){5}}
\put(51,10){\line(0,1){5}}
\put(51,0){\line(0,1){5}}
\put(30,0){\line(1,0){27}}
\put(30,5){\line(1,0){27}}
\put(30,10){\line(1,0){21}}
\put(30,15){\line(1,0){21}}
\put(32,1.5){$t^3_1,t^2_2,t^1_3$}
\put(44,1.5){$t^2_1,t^1_2$}
\put(53,1.5){$t^1_1$}
\put(32,6.5){$t^3_2,t^2_3,t^1_4$}
\put(35.5,11.5){$t^1_6$}
\put(46,11.5){$t^1_5$}
\end{picture}
\end{center}

\vspace{-5mm}

\caption{Example of tableaux for $N=4$ with associated variables $t^a_i$.}\label{fig3}
\end{figure}

In general, for each tableaux $\tchi$ of the weight $\bar n(\tchi)$ and any segment
 $\meg{\bar\lll}{\bar\rr}$, such that $\bar r-\bar l=\bar n$ the set of the variables
$\bar t_{\meg{\bar\lll}{\bar\rr}}$
\begin{gather*}
\bar t_{\seg{\bar\rr}{\bar\lll}} =
\big\{t^{1}_{\lll_{1}+1},\ldots,t^{1}_{\rr_{1}};
t^{2}_{\lll_{2}+1},\ldots,t^{2}_{\rr_{2}};\ldots;
t^{N-1}_{\lll_{N-1}+1},\ldots,t^{N-1}_{\rr_{N-1}} \big\}.
\end{gather*}
decouples into $\heit(\tchi)$ groups of variables
\begin{gather}\label{str-var}
\bar t_{\tchi^i}=
\big\{t^1_{\lll_1+h^{i-1}_1+1},\ldots,t^1_{\lll_1+h^{i}_1}; \ldots ;
t^{N-1}_{\lll_{N-1}+h^{i-1}_{N-1}+1},\ldots,t^{N-1}_{\lll_{N-1}+h^i_{N-1}}\big\}
.
\end{gather}
The variable $t_k^a$ belongs to the subset $\bar t_{\tchi^i}$ if
$\lll_a+h^{i-1}_a<k\leq \lll_a+h^i_a$.
It is located in the $(\lll_a+h^i_a+1-k)$th box of the row $\tchi^i$ counting boxes in this row
from the left edge of the tableaux.

In the same setting we def\/ine
\begin{gather}\label{com-r-s}
Z_{\tchi}(\bar t_{\seg{\bar\rr}{\bar\lll}})=\prod_{1\leq i<j\leq\heit(\tchi)}
Z_{\tchi^i,\tchi^j}(\bar t_{\tchi^i};\bar t_{\tchi^j}),
\end{gather}
where
\begin{gather*}
Z_{\tchi^i,\tchi^j}(\bar t_{\tchi^i};\bar t_{\tchi^j})=
\prod_{a=1}^{N-2}\quad \prod_{\ell=\lll_a+h^{j-1}_a+1}^{\lll_a+h^j_a} \
\prod_{\ell'=\lll_{a+1}+h^{i-1}_{a+1}+1}^{\lll_{a+1}+h^i_{a+1}}
\frac{q-q^{-1} t^a_\ell/t^{a+1}_{\ell'}}{1- t^a_\ell/t^{a+1}_{\ell'}}
\end{gather*}
is a rational series def\/ined by the interchanging of the variables of the type $a+1$
from the $i$-th row and variables of the type $a$ from the $j$th row of the
tableaux $\tchi$.

In our example, the group of variables
\begin{gather}\label{seq1}
\sk{t^{3}_{2},t^3_1;t^2_{3},t^2_2,t^2_1;t^1_{6},t^1_5,t^1_4,t^1_3,t^1_2,t^1_1}
\end{gather}
decomposes into three groups
\begin{gather}\label{seq2}
\sk{\cdot;\cdot;t^1_{6},t^1_5}
\sk{t^{3}_{2};t^2_{3};t^1_4}
\sk{t^3_1;t^2_2,t^2_1;t^1_3,t^1_2,t^1_1} .
\end{gather}
In this example the rational series $Z_{\tchi}(\bar t_{\segg{\bar n}})$ is equal to
\begin{gather}\label{part0}
\prod_{\ell=5,6}\frac{q-q^{-1}\,t^1_\ell/t^{2}_{3}}{1-\,t^1_\ell/t^{2}_{3}}\cdot
\prod_{\ell=4,5,6\atop\ell'=1,2}\frac{q-q^{-1}\,t^1_\ell/t^{2}_{\ell'}}{1- t^1_\ell/t^{2}_{\ell'}}\cdot
\frac{q-q^{-1}\,t^2_3/t^{3}_{1}}{1- t^2_3/t^{3}_{1}} .
\end{gather}

For a given tableaux $\tchi$ we def\/ine  the ordered product
\begin{gather}\label{seq-chi}
\Ccal_{\tchi}(\bar t_{\segg{\bar n}})
=\prod^{\longleftarrow}_{\heit(\tchi)\geq i\geq 1}\Ccal(\bar t_{\tchi^i})
\end{gather}
of the negative projections of the strings \r{coef-m}.
Each factor $\Ccal(\bar t_{\tchi^i})$ corresponds to the  projection of the certain string
depending on the set of variables $\bar t_{\tchi^i}$ \r{str-var}. Decomposition of the tableaux row
into boxes shows the structures of this string. Each box of the type $a$ corresponds to the
composed current $F_{a+1,1}$ which depends on the type 1 variable placed in this box. Variables
of other types  from the same box enter through rational factors.
For the tableaux shown on the Fig.~\ref{fig3} this product reads
\begin{gather*}
\Ccal(t^1_5,t^1_6)\cdot \Ccal(t^1_4;t^2_3;t^3_2)\cdot
\Ccal(t^1_1,t^1_2,t^1_3;t^2_1,t^2_2;t^3_1)
\end{gather*}
and
\begin{gather*}
\Ccal(t^1_5,t^1_6)=\Pfm\sk{F_{2,1}(t^1_6)F_{2,1}(t^1_5)},\qquad
\Ccal(t^1_4;t^2_3;t^3_2)=\frac{1}{1-t^1_4/t^2_3}\frac{1}{1-t^2_3/t^3_2}\Pfm\sk{F_{4,1}(t^1_4)},
\\
\Ccal(t^1_1,t^1_2,t^1_3;t^2_1,t^2_2;t^3_1)=\frac{1}{1-t^1_3/t^2_2}\frac{1}{1-t^2_2/t^3_1}
\frac{1}{1-t^1_2/t^2_1}
\Pfm\sk{F_{4,1}(t^1_3)F_{3,1}(t^1_2)F_{2,1}(t^1_1)}.
\end{gather*}

\begin{proposition}
The coefficient $\Dcal(\bar t_{\segg{\bar n}})$ of the inverse series
$\sk{\Pfm\sk{\tSgs^{N}(u_1,\ldots,u_{N-1})}}^{-1}$ are given
\begin{gather}\label{exp-fact}
\Dcal(\bar t_{\segg{\bar n}})=\sum_{\tchi}(-1)^{\heit(\tchi)+1}
\tSym_{\, \bar t_{\segg{\bar n}}}\sk{
Z_{\tchi}(\bar t_{\segg{\bar n}})\cdot \Ccal_{\tchi}(\bar t_{\segg{\bar n}})}
\end{gather}
by the sum over all possible connected tableaux $\tchi$ such that the weight of tableaux
$\wt(\tchi)$ is equal to $\bar n$.
\end{proposition}

\begin{proof} For arbitrary non-empty set $\bar n\neq\bar 0$
we substitute expression \r{exp-fact} into \r{rec1} to obtain the
relation
\begin{gather}
 \sum_{\tchi\atop \wt(\tchi)=\bar n}(-1)^{\heit(\tchi)+1}
\tSym_{\, \bar t_{\segg{\bar n}}}\sk{
Z_{\tchi}(\bar t_{\segg{\bar n}})\cdot \Ccal_{\tchi}(\bar t_{\segg{\bar n}})}\nonumber\\
\qquad{} +\sum_{\bar m\atop {\bar m}\neq{\bar 0}}
\sum_{\tchi'\atop \wt(\tchi')=\bar n-\bar m}(-1)^{\heit(\tchi')+1}
\tSym_{\, \bar t_{\segg{\bar n}}}\sk{Z_{\bar m}(\bar t_{\segg{\bar n}})\cdot
Z_{\tchi'}(\bar t_{\seg{\bar n}{\bar m}})\cdot  \Ccal_{\tchi'}(\bar t_{\seg{\bar n}{\bar m}})\cdot
\Ccal(\bar t_{\segg{\bar m}})  }\label{pr1}
\end{gather}
which has to be equal 0. After this substitution  the product of the series
$Z_{\tchi'}(\bar t_{\seg{\bar n}{\bar m}})\cdot  \Ccal_{\tchi'}(\bar t_{\seg{\bar n}{\bar m}})$
should be under $q$-symmetrization over the set of the variables
$\bar t_{\seg{\bar n}{\bar m}}$. Since the series
$Z_{\bar m}(\bar t_{\bar n})$ is symmetric with respect to the set of these variables
and the series $\Ccal(\bar t_{\bar m})$ does not depend on the variables
$\bar t_{\seg{\bar n}{\bar m}}$ we may include these
series
under $q$-symmetrization over the variables $\bar t_{\seg{\bar n}{\bar m}}$. Then, since the variables
$\bar t_{\seg{\bar n}{\bar m}}$ forms a subset of the variables
$\bar t_{\segg{\bar n}}$, the $q$-symmetrization over variables
$\bar t_{\seg{\bar n}{\bar m}}$ disappear due to the property \r{sym*}.

We will prove the cancellation of the terms
in \r{pr1} in the sums over tableaux of the f\/ixed height.
Keep the terms in the summation of the f\/irst line of this relation which correspond
to the connected tableaux $\tchi$ such that $n_a(\tchi)=n_a$ and $\heit(\tchi)=p+1$.
Keep the terms in the summation of the second line of \r{pr1} which correspond
to all connected tableaux $\tchi'$ such that $n_a(\tchi')=n_a-m_a$ and $\heit(\tchi')=p$
for f\/ixed $p=0,\ldots,n_1-1$.

Fix a term  in the f\/irst sum of \r{pr1} corresponding to some connected tableaux $\tchi$
with a~weight~$\bar n$.
Consider the f\/irst (bottom) line of the tableaux $\tchi$. Denote by $m_a=c^1_a+\cdots+c^1_{N-1}$
the nonnegative integers def\/ined by this row. It is clear that this set of integers
satisf\/ies the admissibility condition $m_1\geq m_2\geq \cdots \geq m_{N-1}$ and $m_a\leq n_a$.
In the second double sum of~\r{pr1} choose the term corresponding to this set $\bar m$ and
the tableaux  $\tchi'$ def\/ined by the following rule. If we glue from the bottom of the
tableaux $\tchi'$ the row of boxes such that it  has a length $m_1$ and the number of boxes
of the type $a$ is equal to $m_a-m_{a+1}$ then for the obtained
 tableaux $\tilde\chi$ we require $r^i_a(\tilde\chi)=r^i_a(\tchi)$
and $\heit(\tilde\chi)=p+1$ for all possible values $i$ and $a$.
We claim that for each f\/ixed tableaux $\tchi$ there are unique set $\bar m$ such that
$m_a\leq n_a$ and  there are a single tableaux~$\tchi'$ which satisf\/ies above conditions.
The tableaux $\tilde\chi$ and $\tchi$ coincide actually. The product of the coef\/f\/icients
$\Ccal_{\tchi'}(\bar t_{\seg{\bar n}{\bar m}})\cdot
\Ccal(\bar t_{\bar m})$ will be equal to $\Ccal_{\tchi}(\bar t_{\segg{\bar n}})$.

According to the def\/initions of the series
$Z_{\bar m}(\bar t_{\bar n})$ \r{Zserm} and
$Z_{\tchi'}(\bar t_{\seg{\bar n}{\bar m}})$ \r{com-r-s} their product will be equal to the
series $Z_{\tchi}(\bar t_{\segg{\bar n}})$.
The term corresponding to the
f\/ixed tableaux $\tchi$ in the f\/irst line of \r{pr1} and the term from the second line
given by $\bar m$ and $\tchi'$ described above
cancel each other since they will enter with dif\/ferent signs: $\heit(\tchi)=
\heit(\tchi')+1=p+1$. \end{proof}

For the example of the tableaux shown in the Fig.~\ref{fig3} the sets $\bar t_{\segg{\bar m}}$
and  $\bar t_{\seg{\bar n}{\bar m}}$ are
\begin{gather*}
\{t^3_1;t^2_2,t^2_1;t^1_3,t^1_2,t^1_1\}\qquad\mbox{and}\qquad \{t^3_2;t^2_3;t^1_6,t^1_5,t^1_4\}
\end{gather*}
respectively. The series $Z_{\bar m}(\bar t_{\bar n})$ is
\begin{gather}\label{part1}
\prod_{\ell=4,5,6\atop\ell'=1,2}\frac{q-q^{-1} t^1_\ell/t^{2}_{\ell'}}{1- t^1_\ell/t^{2}_{\ell'}}\cdot
\frac{q-q^{-1} t^2_3/t^{3}_{1}}{1- t^2_3/t^{3}_{1}} .
\end{gather}
The series $Z_{\tchi'}(\bar t_{\seg{\bar n}{\bar m}})$ is
\begin{gather}\label{part2}
\prod_{\ell=5,6}\frac{q-q^{-1} t^1_\ell/t^{2}_{3}}{1- t^1_\ell/t^{2}_{3}} .
\end{gather}
The product of \r{part1} and \r{part2} obviously coincides with~\r{part0}.

\subsection[Inversion of generating series for $U_q(\widehat{\mathfrak{gl}}_2)$]{Inversion of generating series for $\boldsymbol{\Uqgl{2}}$}

Quantum af\/f\/ine algebra $\Uqgl{2}$ in its current realization formed by the modes
of the currents\footnote{Since algebra $\mathfrak{gl}_2$ has only one root
we remove index of this simple root in the notation of the currents in
case of the current realization of the algebra $\Uqgl{2}$.} $E(t)$, $F(t)$ and
Cartan currents $k_1^\pm(t)$, $k^\pm_2(t)$.

Let us invert explicitly the generating series $\Pfm\sk{\Fgs(\bar u)}$ in the simplest case
of one generating parameters $\bar u=u_1\equiv u$ in \r{n-o-gs} to show what kind of relations
can be obtained for the generating series of the $\Uqgl{2}$
of\/f-shell Bethe vectors $\Pfp\sk{\Fgs(u)}$:
\begin{gather}\label{sp-gs}
\Pfp\sk{\Fgs(u)}=\Pfm\sk{\Fgs(u)}^{-1} \star \Fgs(u) .
\end{gather}

For any positive $n$  and non-negative $p\leq n$ we def\/ine the set of $p+1$ positive integers
$\{\bar k_p\}=\{k_1,k_2,\ldots,k_p,k_{p+1}\}$ such that $0=k_0<k_1<k_2<\cdots<k_p<k_{p+1}=n$.
Using this data we def\/ine the ordered products
\begin{gather*}
\bFF_{\bar k_p}(\bar t_{\segg{n}})=
\prod^{\longleftarrow}_{p+1\geq m\geq 1}
\frac{1}{(k_m-k_{m-1})!}  \Pfm\sk{F(t_{k_m})\cdots F(t_{k_{m-1}+1}) } .
\end{gather*}
It is clear that the inverse generating series $\Pfm\sk{\Fgs(u)}$ can be written using
$q$-symmetrization of these
ordered products $\bFF_{\bar k_p}(\bar t_{\segg{n}})$ as follows
\begin{gather}\label{nigs}
\Pfm\sk{\Fgs(u)}^{-1}=1+\sum_{n> 0}  \tSym_{\, \bar t_{\segg{n}}}\sk{
\sum_{p=0}^{n-1}(-1)^{p+1}\sum_{\{\bar k_p\}}
\bFF_{\bar k_p}(\bar t_{\segg{n}})
} u^n .
\end{gather}

Then using the def\/inition of the $\star$-product we can obtain from \r{sp-gs}
a special presentation for the $\Uqgl{2}$ universal weight function
\begin{gather}\label{DF-gen}
\Pfp\sk{\Fgs(\bar t_{\segg{n}})}=\tSym_{\, \bar t_{\segg{n}}}\sk{
\sum_{s=0}^n\sum_{p=0}^{n-s-1}(-1)^{p+1}\sum_{\{k_p\}}
  \bFF_{\bar k_p}(\bar t_{\seg{n}{s}})\cdot \Fgs(\bar t_{\segg{s}})  },
\end{gather}
where summation over the set $\{\bar k_p\}$ runs over all possible $k_i$ such that
$s=k_0<k_1<k_2<\cdots<k_p<k_{p+1}=n$. An extreme term in the sum when $s=n$  and the sum
over $p$ is absent corresponds to the product of the total currents $F(t_n)\cdots F(t_1)$.
Note that an equality~\r{DF-gen} can be treated as generalization of Ding--Frenkel
relation $\Pfp\sk{F(t)}=F(t)-\Pfm\sk{F(t)}$ \cite{DF} when the positive projection is taken
from the product of the currents.

\subsection[Universal Bethe ansatz for $U_q(\widehat{\mathfrak{gl}}_2)$]{Universal Bethe ansatz for $\boldsymbol{\Uqgl{2}}$}\label{gl2-ex}

Let $J$ be the left ideal of $\Uqbp\subset\Uqgl{2}$, generated by all element of the form
$\Uqbp\cdot E[n]$,  $n>0$.
As it was mentioned above the standard Borel subalgebra $\Uqbp$ in terms of the currents generators
is formed by the modes $F[n]$, $k^+_{1,2}[n]$, $n\geq 0$ and $E[m]$, $m>0$. Let us denote subalgebras
generated by these modes as $U_f^+$, $U^+_k$ and $U^+_e$, respectively. The multiplication in~$\Uqbp$ implies  an isomorphism
of the vectors spaces
\begin{gather*}
U_f^+\otimes U^+_k\otimes U^+_e\to \Uqbp .
\end{gather*}

We introduce ordering of the generators in the Borel subalgebra $\Uqbp$
\begin{gather}\label{ordering}
U_f^+ \prec  U^+_k\prec U^+_e
\end{gather}
induced by the circular ordering of the Cartan--Weyl generators in the whole algebra
$\Uqgl{2}$ \cite{EKhP}.
We call any element $w\in\Uqbp$ {\it normal ordered} and denote it as $:W:$ if it is presented
as the linear combination of the elements of the form $W_1\cdot W_2\cdot W_3$, where
$W_1\in U_f^+$, $W_2\in U^+_k$, $W_3\in U^+_e$. It is convenient to gather the generators of the
subalgebras $U^+_f$ and $U^+_e$ into generating series
\begin{gather*}
F^+(t)=\sum_{n\geq 0}F[n]t^{-n},\qquad E^+(t)=\sum_{n> 0}E[n]t^{-n}
\end{gather*}
which we call the half-currents.

{\em A universal transfer matrix} is the following combinations of the Cartan and half-currents
\begin{gather}\label{utm}
\mathcal{T}(t)=k_1^+(t)+ F^+(t)  k^+_2(t)  E^+(t) +k^+_2(t).
\end{gather}
Using the commutation relations in the algebra $\Uqbp\subset\Uqgl{2}$ one may check that
these transfer matrices commute for the dif\/ferent values of the spectral parameters
\begin{gather*}
[\mathcal{T}(t),\mathcal{T}(t')]=0
\end{gather*}
and so generates the inf\/inite set of commuting quantities\footnote{A standard
way to prove this commutativity is to note that \r{utm} is a trace of the fundamental
$\LL$-operators for $\Uqgl{2}$ and the commutativity follows from the Yang--Baxter
equation for these $\LL$-operators.}. We are interesting in the ordering
relations
between universal transfer matrix $\mathcal{T}(t)$ and the universal Bethe vector
$\Pfp\sk{\Fgs(\bar t_{\segg{n}})}$.
Note that the universal transfer matrix is ordered according to the ordering~\r{ordering}.

\begin{proposition}\label{main-th}
A formal series identity is valid in $\Uqbp$
\begin{gather}\label{main-id}
:\mathcal{T}(t)\cdot \Pfp\sk{\Fgs(\bar t_{\segg{n}})}:\ \,
=\Pfp\sk{\Fgs(\bar t_{\segg{n}})}\cdot\tau(t;\bar{t}_{\segg{n}})
 \quad{\rm mod}\ J
\end{gather}
modulo elements of the left ideal $J$
if the set $\{t_j\}$ of the Bethe parameters satisfies the set of the universal Bethe equations~{\rm \cite{ACDFR}},
$j=1,\ldots,n$:
\begin{gather}\label{univ-BE}
\frac{k^+_1(t_j)}{k^+_{2}(t_j)}=
\prod_{m\neq j}^{n}\frac{qt_j-q^{-1}t_m}{q^{-1}t_j-qt_m}
\end{gather}
and
\begin{gather}\label{univ-ev}
\tau(t;\bar{t}_{\segg{n}})
= k^+_1(t) \prod_{j=1}^{n}\frac{q^{-1}t-qt_j}{t-t_j}
+ k^+_2(t) \prod_{j=1}^{n}\frac{qt-q^{-1}t_j}{t-t_j}
\end{gather}
is an eigenvalue of the universal transfer matrix.
\end{proposition}

\begin{proof} Recall that a universal Bethe vector in the considered case
coincides with projection of the product of the currents:
$\Pfp\sk{\Fgs(\bar t_{\segg{n}})}=\Pfp\sk{F(t_n)\cdots F(t_1)}/n!$ and can be
presented as factorized product of the linear combinations of the half-currents
$F^+(t_i)$ \cite{EKhP,KP}.
 A direct way to prove the statement of the Proposition~\ref{main-th} is to use
the commutation relations in $\Uqbp$ between half-currents
\begin{gather*}
[E^+(t),F^+(t')]= \frac{(q-q^{-1})t'}{t-t'}  \sk{k_1^+(t')k_2^+(t')^{-1}- k_1^+(t)k_2^+(t)^{-1}} ,
\\
k_2^+(t) F^+(t') k_2^+(t)^{-1}= \frac{qt-q^{-1}t'}{t-t'} F^+(t')-
\frac{(q-q^{-1})t'}{t-t'} F^+(t)
\end{gather*}
and similar for $k_1^+(t) F^+(t') k_1^+(t)^{-1}$
to present the product $\mathcal{T}(t)\cdot \Pfp\sk{\Fgs(\bar t_{\segg{n}})}$
in the normal ordered form. But this way is not easy even in the simplest case
of the algebra $\Uqgl{2}$. It becomes much more involved in the general case
of the algebra $\Uqgln$. There is a simple way to avoid these dif\/f\/iculties using
the relation~\r{DF-gen}.

This relation allows to replace the  projection of the product of the currents onto
positive Borel subalgebra $\Uqbp$ by the linear combination of the terms
\begin{gather}
\Pfp\sk{F(t_{n})\cdots F(t_{1})}= F(t_{n})\cdots F(t_{1})\nonumber\\
\phantom{\Pfp\sk{F(t_{n})\cdots F(t_{1})}=} {}- n  \tSym_{\, \bar t}
\sk {\Pfm\sk{F(t_{n})}\cdot F(t_{n-1})\cdots F(t_{1})} + W ,\label{sef2}
\end{gather}
where $W$ are the terms which have on the left the product of at least  two negative projections
of the currents $F(t)$. Let us substitute \r{sef2} into the product
$\mathcal{T}(t)\cdot \Pfp\sk{\Fgs(\bar t_{\segg{n}})}$. Although left hand side of \r{sef2}
belongs to the positive Borel subalgebra, each term in the right hand side of the equality
\r{sef2} does not belong to $\Uqbp$. The product of the universal transfer matrix
with these terms will produce the terms such that some of them belong to $\Uqbp$ and other does not
belong to $\Uqbp$. The latter
terms which after ordering do not belong to $\Uqbp$  can be omitted since we are interesting
only in the terms which belong to the positive Borel subalgebras. In particular, one may
check that $:\mathcal{T}(t)\cdot W:\,\not\in\Uqbp$, where $W$ are the terms in \r{sef2}
not showing explicitly (see \cite{FKPR} for details).
Using the commutation relations of the half-currents $E^+(t)$ with total currents $F(t')$
\begin{gather*}
[E^+(t),F(t')]= \frac{(q-q^{-1})t'}{t-t'}  \sk{k_1^+(t')k_2^+(t')^{-1}- k_1^-(t')k_2^-(t')^{-1}} ,
\end{gather*}
the commutation relations of the Cartan currents with total currents $F(t')$
\begin{gather*}
k_2^+(t) F(t') k_2^+(t)^{-1}= \frac{qt-q^{-1}t'}{t-t'} F^+(t') ,\qquad
k_1^+(t) F(t') k_1^+(t)^{-1}= \frac{q^{-1}t-qt'}{t-t'} F^+(t')
\end{gather*}
and the commutation relations of $E^+(t)$, $k^+_1(t)$, $k^+_2(t)$ with
$\Pfm\sk{F(t')}$
we may check that the only terms which belong to the positive Borel subalgebra $\Uqbp$
and do not belong to the left ideal $J$  in the normal ordered product
$:\mathcal{T}(t)\cdot \Pfp\sk{F(t_{n})\cdots F(t_{1})}:$ are
\begin{gather}
 :\mathcal{T}(t)\cdot \Pfp\sk{F(t_{n})\cdots F(t_{1})}:\,
 \nonumber\\
 \quad= \Pfp\sk{F(t_{n})\cdots F(t_{1})}\sk{\prod_{i=1}^n\frac{q^{-1}t-qt_i}{t-t_i}  k^+_1(t)+
\prod_{i=1}^n\frac{qt-q^{-1}t_i}{t-t_i}  k^+_2(t)} \nonumber\\
 \qquad{}+ n  \overline{\rm Sym}_{\, \bar t}
\sk{\Pfp\sk{F^+(t)k^+_2(t)  F(t_{n})\cdots F(t_{2})}
\frac{(q-q^{-1})t_1}{t-t_1}   k_1^+(t_1)k_2^+(t_1)^{-1}  } \nonumber\\
 \qquad{}- n  \overline{\rm Sym}_{\, \bar t}
\sk{\Pfp\sk{F^+(t)k^+_2(t)  F(t_{n-1})\cdots F(t_{1})}
\frac{(q-q^{-1})t_n}{t-t_n}  } .\label{lhs3}
\end{gather}
In order to prove the statement of the Proposition~\ref{main-th} we have to cancel
the last two terms in~\r{lhs3}. This can be done using the properties of the $q$-symmetrization
that for any formal series $G(t_1,\ldots,t_n)$ on $n$ formal
variables $t_i$ we have
\begin{gather*}\label{sym-last}
n  \overline{\rm Sym}_{\, \bar t}\ G(t_1,\ldots,t_n)=
\sum_{m=1}^n \prod_{j=m+1}^n \frac{q^{-1}t_m-qt_j}{qt_m-q^{-1}t_j} \,
\overline{\rm Sym}_{\, \bar t\setminus t_m}  G(t_1,\ldots,t_{m-1},t_{m+1},\ldots,t_n,t_m)
\end{gather*}
and
\begin{gather*}\label{sym-first}
n  \overline{\rm Sym}_{\, \bar t}   G(t_1,\ldots,t_n)=
\sum_{m=1}^n \prod_{j=1}^{m-1} \frac{q^{-1}t_j-qt_m}{qt_j-q^{-1}t_m}
\overline{\rm Sym}_{\, \bar t\setminus t_m}\ G(t_m,t_1,\ldots,t_{m-1},t_{m+1},\ldots,t_n) ,
\end{gather*}
where $q$-symmetrization in the right hand sides of this formal series identities
runs over $(n-1)$ variables $\bar t\setminus t_m=\{t_1,\ldots,t_{m-1},t_{m+1},\ldots,t_m\}$.
Using these relations we conclude that last two terms in \r{lhs3} cancel each other
provided the relation \r{univ-BE} is satisf\/ied.\end{proof}

\appendix

\section[Current realization of $U_q(\widehat{\mathfrak{gl}}_N)$]{Current realization of $\boldsymbol{\Uqgln}$}\label{com-rel}

The commutation relations for the algebra $\Uqgln$ in the current realization are
given by the following set of the relations
\begin{gather}
(q^{-1}z-q^{}w)E_{i}(z)E_{i}(w)= E_{i}(w)E_{i}(z)(q^{}z-q^{-1}w)  ,\nonumber
\\
(z-w)E_{i}(z)E_{i+1}(w)= E_{i+1}(w)E_{i}(z)(q^{-1}z-qw)  ,\nonumber
\\
k_i^\pm(z)E_i(w)\left(k_i^\pm(z)\right)^{-1}=
\frac{z-w}{q^{-1}z-q^{}w}E_i(w)  ,\nonumber
\\
k_{i+1}^\pm(z)E_i(w)\left(k_{i+1}^\pm(z)\right)^{-1}=
\frac{z-w}{q^{}z-q^{-1}w}E_i(w)  ,\nonumber
\\
k_i^\pm(z)E_j(w)\left(k_i^\pm(z)\right)^{-1}=E_j(w),
\qquad {\rm if}\quad i\not=j,j+1  ,\nonumber
\\
 (q^{}z-q^{-1}w)F_{i}(z)F_{i}(w)= F_{i}(w)F_{i}(z)(q^{-1}z-q^{}w)
 ,  \label{gln-com1}
\\
(q^{-1}z-qw)F_{i}(z)F_{i+1}(w)= F_{i+1}(w)F_{i}(z)(z-w)  ,\nonumber
\\
k_i^\pm(z)F_i(w)\left(k_i^\pm(z)\right)^{-1}=
\frac{q^{-1}z-qw}{z-w}F_i(w)  ,\nonumber
\\
k_{i+1}^\pm(z)F_i(w)\left(k_{i+1}^\pm(z)\right)^{-1}=
\frac{q^{}z-q^{-1}w}{z-w}F_i(w)  ,\nonumber
\\
k_i^\pm(z)F_j(w)\left(k_i^\pm(z)\right)^{-1}=F_j(w),
 \qquad {\rm if}\quad i\not=j,j+1 ,\nonumber
\\
[E_{i}(z),F_{j}(w)]= \delta_{{i},{j}}\ \delta(z/w)\
(q-q^{-1})\left( k^+_{i}(z)/k^+_{i+1}(z)-k^-_{i}(w)/k^-_{i+1}(w)\right)\nonumber
\end{gather}
and the Serre relations for the  currents $E_{i}(z)$ and $F_{i}(z)$
\begin{gather}
{\rm Sym}_{z_1,z_{2}}
(E_{i}(z_1)E_{i}(z_2)E_{i\pm 1}(w)
 -(q+q^{-1})E_{i}(z_1)E_{i\pm 1}(w)E_{i}(z_2) \nonumber\\
\qquad{}+E_{i\pm 1}(w)E_{i}(z_1)E_{i}(z_2))=0 ,\nonumber\\
{\rm Sym}_{z_1,z_{2}}
(F_{i}(z_1)F_{i}(z_2)F_{i\pm 1}(w)
 -(q+q^{-1})F_{i}(z_1)F_{i\pm 1}(w)F_{i}(z_2)\nonumber\\
\qquad{}+F_{i\pm 1}(w)F_{i}(z_1)F_{i}(z_2))=0  .\label{serre}
\end{gather}
Formulae \eqref{gln-com1} and \eqref{serre} should be considered as
formal series identities describing the inf\/inite set of the relations between modes of the
currents. The symbol $\delta(z)$ entering these relations is a formal series~$\sum_{n\in\ZZ} z^n$.

Following \cite{DKh,KP}, we introduce  {\it composed currents\/}
$\ff_{j,i}(t)$ for $i<j$. The composed currents for
 nontwisted quantum af\/f\/ine algebras were def\/ined in~\cite{DKh}.
According to this paper, the coef\/f\/icients of the series $\ff_{j,i}(t)$ belong to
the completion $\overline U_F$ of the algebra $U_F$.

The completion $\overline U_F$ determines analyticity properties of products
of currents (and coincide with analytical properties of their matrix coef\/f\/icients for highest
weight representations \cite{E}). One can show that for $|i-j|>1$, the product
$F_i(t)F_j(w)$ is an expansion of a function analytic at $t\ne 0$, $w\ne 0$.
The situation is more delicate for $j=i,i\pm1$. The products $F_i(t)F_i(w)$ and
$F_i(t)F_{i+1}(w)$ are expansions of analytic functions at $|w|<|q^2 t|$, while
the product $F_i(t)F_{i-1}(w)$ is an expansion of an analytic function at
$|w|<|t|$. Moreover, the only singularity of the corresponding functions
in the whole region $t\ne 0$, $w\ne 0$, are simple poles at the respective
hyperplanes, $w=q^2t$ for $j=i,i+1$, and $w=t$ for $j=i-1$. Recall, that the deformation
parameter $q$ is a generic complex number, which is neither 0 nor a root of unity.

The def\/inition of the composed currents may be written in analytical form
\begin{gather}\label{rec-f1}
\ff_{j,i}(t)  =
- \mathop{\rm res}\limits_{w=t}\ff_{j,a}(t)\ff_{a,i}(w) \frac{dw}{w}  =
\mathop{\rm res}\limits_{w=t}\ff_{j,a}(w)\ff_{a,i}(t) \frac{dw}{w}
\end{gather}
for any $a=i+1,\ldots, j-1$. It is equivalent to the relation
\begin{gather}
\ff_{j,i}(t) =
\oint \ff_{j,a}(t) \ff_{a,i}(w)  \frac{dw}{w}-
\oint \frac{q^{-1}-qt/w}{1-t/w}
 \ff_{a,i}(w) \ff_{j,a}(t)  \frac{dw}{w} ,\nonumber\\
\ff_{j,i}(t) =
\oint \ff_{j,a}(w) \ff_{a,i}(t)  \frac{dw}{w}-
\oint  \frac{q^{-1}-qw/t}{1-w/t}
 \ff_{a,i}(t) \ff_{j,a}(w)  \frac{dw}{w} . \label{rec-f111}
\end{gather}
In \r{rec-f111} $\oint \frac{dw}{w} g(w)=g_0$ for any formal series $g(w)=\sum_{n\in\ZZ}g_n z^{-n}$.

Using the  relations \r{gln-com1} on $F_i(t)$ we can calculate
the residues in \r{rec-f1} and  obtain the following expressions  for  $\ff_{j,i}(t)$, $i<j$:
\begin{gather}\label{res-in}
\ff_{j,i}(t)=(q-q^{-1})^{j-i-1}\ff_i(t) \ff_{i+1}(t)\cdots \ff_{j-1}(t) .
\end{gather}
For example, $F_{i+1,i}(t)=F_i(t)$, and
$F_{i+2,i}(t)=(q-q^{-1})F_i(t)F_{i+1}(t)$. The last product is well-def\/ined
according to the analyticity properties of the product $F_i(t)F_{i+1}(w)$,
described above. In a similar way, one can show inductively that the product
in the right hand side of \r{res-in} makes sense for any $i<j$.
Formulas \r{res-in} prove that the def\/ining relations for the composed currents~\r{rec-f1} or \r{rec-f111}  yields the same answers for all possible
values $i<a<j$.

Calculating formal integrals in \r{rec-f111} we obtain
the following presentations
 for the composed currents:
\begin{gather*}\label{rec-f112a}
F_{j,i}(t)=
F_{j,a}(t)F_{a,i}[0]-q^{-1}F_{a,i}[0]F_{j,a}(t)+
(q-q^{-1})\sum_{k< 0}F_{a,i}[k]\,F_{j,a}(t)\,t^{-k} ,
\\
\label{rec-f112b}
F_{j,i}(t)=
F_{j,a}[0]F_{a,i}(t)-qF_{a,i}(t)F_{j,a}[0]
+(q-q^{-1})\sum_{k\geq 0}F_{a,i}(t)\,F_{j,a}[k]\,t^{-k},
\end{gather*}
which are useful for the calculation of their projections.

\subsection*{Acknowledgements}

The main idea to use generating series for the description of the hierarchical
Bethe ansatz appeared during authors visit to Max-Planck Institute f\"ur Mathematik,
Bonn, in January, 2008. Authors acknowledge this scientif\/ic center for the hospitality and
stimulating scientif\/ic atmosphere.

This work was partially done when the second
author (S.P.) visited Laboratoire d'An\-necy-Le-Vieux de Physique Th\'eorique in 2006 and 2007.
These visits were possible due to the f\/inancial support of
the CNRS-Russia exchange program on mathematical physics.
He thanks LAPTH for the hospitality and stimulating scientif\/ic atmosphere.
Authors are grateful to Luc Frappat and \'Eric Ragoucy for many helpful discussions.
The authors were supported in part by RFBR
grant  08-01-00392 and grant for the support
of scientif\/ic schools NSh-3036.2008.2. The f\/irst author was also supported
 by the Atomic Energy Agency of the Russian Federation,
and by the ANR grant 05-BLAN-0029-01.
The second author  was also supported in part by RFBR-CNRS grant 07-02-92166-CNRS.

\pdfbookmark[1]{References}{ref}
\LastPageEnding


\begin{thebibliography}{99}

\footnotesize\itemsep=0pt



\bibitem{ACDFR}
Arnaudon D.,  Cramp\`e N., Doikou A.,  Frappat L., Ragoucy E.,
Spectrum and Bethe ansatz equations for the
$U_{q} {\left( {gl(\mathcal{N})} \right)}$ closed and open spin chains in any representation,
{\it  Ann. Henri  Poincar\'e} {\bf 7} (2006), 1217--1268, \href{http://arxiv.org/abs/math-ph/0512037}{math-ph/0512037}.

\bibitem{D88} Drinfel'd V.G., New realization of Yangians and quantum
af\/f\/ine algebras, {\it Soviet Math. Dokl.} {\bf 36} (1988), 212--216.

\bibitem{DF} Ding J.T., Frenkel I.B., Isomorphism of two realizations of
quantum af\/f\/ine algebra $U_q(\widehat{\mathfrak{gl}}_N)$, {\it Comm. Math. Phys.} \textbf{156} (1993),
277--300.


\bibitem{DKh} Ding J., Khoroshkin S., Weyl group extension of quantized
current algebras, {\it Transform. Groups} {\bf 5} (2000), 35--59, \href{http://arxiv.org/abs/math.QA/9804139}{math.QA/9804139}.


\bibitem{E} Enriquez B., On correlation functions of Drinfeld
currents and shuf\/f\/le algebras, {\it Transform. Groups} {\bf 5} (2000),
111--120, \href{http://arxiv.org/abs/math.QA/9809036}{math.QA/9809036}.

\bibitem{EKhP} Enriquez B., Khoroshkin S., Pakuliak S., Weight
functions and Drinfeld currents, {\it Comm. Math. Phys.} {\bf 276} (2007), 691--725, \href{http://arxiv.org/abs/math.QA/0610398}{math.QA/0610398}.

\bibitem{ER} Enriquez B., Rubtsov V., Quasi-Hopf algebras associated with
$\mathfrak{sl}_2$ and complex curves, {\it Israel J. Math.} {\bf 112} (1999),
61--108, \href{http://arxiv.org/abs/q-alg/9608005}{q-alg/9608005}.

\bibitem{FKPR} Frappat L.,  Khoroshkin S.,  Pakuliak S., Ragoucy \'E.,
Bethe ansatz for the universal weight function,
\href{http://arxiv.org/abs/0810.3135}{arXiv:0810.3135}.


\bibitem{KP} Khoroshkin S., Pakuliak S., The weight function for the quantum af\/f\/ine algebra
$U_q(\widehat{\mathfrak{sl}}_3)$, {\it Theor. and
Math. Phys.} {\bf 145} (2005), 1373--1399, \href{http://arxiv.org/abs/math.QA/0610433}{math.QA/0610433}.

\bibitem{KPT} Khoroshkin S., Pakuliak S., Tarasov V., Of\/f-shell Bethe vectors
and Drinfeld currents, {\it J. Geom. Phys. } {\bf 57} (2007), 1713--1732, \href{http://arxiv.org/abs/math.QA/0610517}{math.QA/0610517}.

\bibitem{KhP-GLN}  Khoroshkin S., Pakuliak S.,
    A computation of an universal weight function for the quantum af\/f\/ine algebra $U_q(\widehat{gl}_{N})$,
    {\it J. Math. Kyoto Univ.} {\bf 48} (2008), 277--322,
\href{http://arxiv.org/abs/0711.2819}{arXiv:0711.2819}.


\bibitem{KR83} Kulish P., Reshetikhin N., Diagonalization of $GL(N)$
invariant transfer matrices and quantum $N$-wave system (Lee model),
{\it J.~Phys.~A: Math. Gen.} {\bf 16} (1983), L591--L596.

\bibitem{MTV} Mukhin E., Tarasov V., Varchenko A.,
Bethe eigenvectors of higher transfer matrices,
{\it J. Stat. Mech. Theory Exp.} {\bf 2006} (2006), no.~8, P08002, 44~pages, \href{http://arxiv.org/abs/math.QA/0605015}{math.QA/0605015}.


\bibitem{OPS} Oskin A., Pakuliak S., Silantyev A.,
On the universal weight function for the quantum af\/f\/ine
    algebra $U_q(\widehat{\mathfrak{gl}}_{N})$,
\href{http://arxiv.org/abs/0711.2821}{arXiv:0711.2821}.


\bibitem{RS} Reshetikhin N., Semenov-Tian-Shansky M., Central extentions
of quantum current groups, {\it Lett. Math. Phys.} {\bf 19} (1990), 133--142.

\bibitem{VT} Tarasov V., Varchenko A., Jackson integrals for the
solutions to Knizhnik--Za\-mo\-lod\-chi\-kov equation, {\it St. Petersburg Math. J.} {\bf
2} (1995), no.~2, 275--313.

\bibitem{VT1} Tarasov V., Varchenko A.,
Geometry of $q$-hypergeometric functions, quantum af\/f\/ine algebras
and elliptic quantum groups, {\it Ast\'erisque} {\bf 246} (1997), 1--135, \href{http://arxiv.org/abs/q-alg/9703044}{q-alg/9703044}.



\bibitem{VT2} Tarasov V., Varchenko A.,
Combinatorial formulae for nested Bethe vectors,
\href{http://arxiv.org/abs/math.QA/0702277}{math.QA/0702277}.

\end{thebibliography}
\end{document}